\documentclass[11pt,reqno,final]{amsart}

\usepackage[utf8]{inputenc}
\usepackage[theoremdefs]{latexdev}
\usepackage{tcolorbox}
\usepackage{todonotes}
\usepackage{algorithmic}
\usepackage{algorithm}
\usepackage{xcolor}
\usepackage{mdframed}

\usepackage{geometry}                
\usepackage[numbers]{natbib}

\DeclareMathOperator{\grad}{grad}

\newcommand{\pair}[1]{\left\langle #1 \right\rangle}

\providecommand{\abs}[1]{\lvert#1\rvert}
\providecommand{\vect}[1]{\boldsymbol{#1}}

\newcommand{\ud}{\mathrm{d}}

\newcommand{\RR}{{\mathbb R}}

\newcommand{\vol}{\ud x}

\newcommand{\Diff}{\mathrm{Diff}}
\newcommand{\Xcal}{\mathfrak{X}}

\newcommand{\SDiff}{{\mathrm{SDiff}}}

\newcommand{\Dens}{\mathrm{Dens}}

\newcommand{\Mass}{\mathrm{Mass}}

\newcommand*\id{\mathrm{id}}

\newcommand{\Densalpha}{\operatorname{Dens}^\alpha}
\newcommand{\GI}{G^I}
\newcommand{\dd}[2]{\frac{d#1}{d#2}}

\DeclareMathOperator*{\argmin}{arg\,min}

\title[Diffeomorphic density registration]{Diffeomorphic density registration}

\author[Bauer]{Martin Bauer}
\address[M.\ Bauer]{Department of Mathematics, Florida State University}
\email{bauer@math.fsu.edu}

\author[Joshi]{Sarang Joshi}
\address[S. Joshi]{Department of Bioengineering, Scientific Computing and Imaging Institute, University of Utah}
\email{scjoshi@gmail.com}

\author[Modin]{Klas Modin}
\address[K.\ Modin]{Department of Mathematical Sciences, Chalmers University of Technology and the University of Gothenburg}
\email{klas.modin@chalmers.se}

\date{\today}                                           

\keywords{density registration, information geometry, Fisher-Rao metric, optimal transport, image registration, diffeomorphism groups, random sampling
}
\subjclass[2010]{58B20, 58E10, 35G25, 35Q31, 76N10}
\begin{document}

\begin{abstract}In this book chapter we study the Riemannian Geometry of the density registration problem: Given two densities  (not necessarily probability densities) defined on a smooth finite dimensional manifold find a diffeomorphism which transforms one to the other. This problem is motivated by the medical imaging application of tracking organ motion due to respiration in Thoracic CT imaging where the fundamental physical property of conservation of mass naturally leads to modeling CT attenuation as a density.

We will study the intimate link between the Riemannian metrics on the space of diffeomorphisms and those on the space of densities. 
We finally develop novel computationally efficient algorithms and demonstrate there applicability for registering RCCT thoracic imaging.
\end{abstract}

\maketitle

\tableofcontents


\section{Introduction}
 Over the last decade image registration has received intense interest, both with respect to medical imaging applications as well to the mathematical foundations of the general problem of estimating a transformation that brings two or more given medical images into a common coordinate system~\cite{grenander1998computational,joshi2004unbiased,younes2010shapes,miller2002metrics,vercauteren2009diffeomorphic,avants2011reproducible,ashburner2007fast,davis2010population}. In this chapter we focus on a subclass of registration problems which term as density registration. The primary difference between density registration and general image registration is in how the registration transformation acts on the image being transformed. In density registration the transformation not only deforms the underlying coordinate system but also scales the image intensity by the local change in volume to conserve ``mass''. In numerous medical imaging applications this conservation of mass is of critical importance and is a fundamental property of the registration problem. The primary motivating clinical application is that of estimating the complex changes in anatomy due to breathing as imaged via 4D respiratory correlated computed tomography (4DRCCT). When tracking organ motion due to breathing it is only natural to assume that the overall mass of the subject being imaged is conserved. Given the physical quantitative nature of CT imaging the natural action of a transformation on a CT image is that of density action: Any local compression induces a corresponding change in local density resulting in changes in the local attenuation coefficient. We will also see that this difference in action of the transformation on the image being registered has wide ranging implications to the structure of the estimation problem. In this chapter we will study both the fundamental geometrical structure of the problem as well as exemplify it's application. The basic outline is as follows: We will first study the abstract mathematical structure of the problem, precisely defining the space of densities and the space of transformation. We will also study the set of transformations that leave the density unchanged. We will see that the explicit characterization of this set of transformations plays a critical role in understanding the geometric structure of the density registration problem.
We will then introduce the general (regularized) density matching problem and present efficient numerical algorithms for several specific choices of regularizers. Finally we will present the before mentioned application to model breathing as imaged via 4D respiratory correlated computed tomography.

\section{Diffeomorphisms and densities}
Let $M$ denote a smooth oriented Riemannian manifold of dimension~$n$ with (reference) volume form $\vol$.

\begin{definition}\label{def:dens}
	The space of \emph{smooth densities}\footnote{If $M$ is compact, then $\Dens(M)$ is an infinite-dimensional \emph{Fr\'echet manifold}~\cite{Ha1982}, i.e., a manifold modelled on a Fr\'echet space. For the purpose of analysis, it is often useful to instead work with the Sobolev completion $\Dens^s(M)$ for a Sobolev index $s>n/2$. $\Dens^s(M)$ is then a \emph{Banach manifold}~\cite{La1999}. The benefit of Banach over Fr\'echet manifolds is that most standard results from finite dimensions, such as the inverse function theorem, are valid.}
	on $M$ is given by
	\begin{equation}
		\Dens(M) = \{ \rho \in C^\infty(M)\mid \rho(x)>0 \quad \forall\, x\in M\}.
	\end{equation}
	The \emph{mass} of a subset $\Omega\subset M$ with respect to $\rho\in\Dens(M)$ is given by $$\Mass_\rho(\Omega) = \int_\Omega \rho\,\vol.$$
\end{definition}



As the focus of the chapter is the registration of densities via transformation, the group $\Diff(M)$ of smooth diffeomorphisms of the manifold plays a central role. 

\begin{definition}\label{def:diff}
	The set of diffeomorphisms on $M$, denoted $\Diff(M)$, consists of smooth bijective mappings $M \to M$ with smooth inverses.
	This set has a natural group structure by composition of maps.
	The Lie algebra of $\Diff(M)$ is given by the space~$\Xcal(M)$ of smooth vector fields (tangential if $M$ has a boundary).\footnote{If $M$ is compact, then $\Diff(M)$ is a \emph{Fr\'echet Lie group}~\cite[\S\!~I.4.6]{Ha1982}, i.e., a  Fr\'echet manifold where the group operations are smooth mappings.}
\end{definition}

The group of diffeomorphisms acts naturally on the space of densities via pullback and pushforward of densities.
Indeed, \emph{pullback of densities} is given by
\begin{equation}
	\Diff(M)\times\Dens(M) \ni (\varphi,\rho) \mapsto \varphi^*\rho = \abs{D\varphi}\,\rho(\varphi(\cdot)) .
\end{equation}
Notice that this is a right action, i.e., $(\varphi\circ\eta)^*\rho = \eta^*(\varphi^*\rho)$.
The corresponding left action is given by \emph{pushforward of densities}
\begin{eqnarray*}
	\Diff(M)\times\Dens(M) \ni (\varphi,\rho) \mapsto \varphi_*\rho  &=& (\varphi^{-1})^*\rho \\
	 &=& \abs{D\varphi^{-1}}\, \rho(\varphi^{-1}(\cdot)).
\end{eqnarray*}

The action of $\Diff(M)$ on densities captures the notion of conservation of mass and is fundamentally different from the standard action of $\Diff(M)$ on functions given by composition.
Indeed, for the density action we have, for any subset $\Omega\subset M$
\begin{equation}
 	\Mass_{\varphi_*\rho}(\varphi(\Omega)) = \Mass_{\rho}(\Omega),
\end{equation}
which follows from the change of coordinates formula for integrals.

The \emph{isotropy subgroup} of an element $\rho\in\Dens(M)$ is, by definition, the subgroup of $\Diff(M)$ which leaves the density $\rho$ unchanged.
It is given by
\begin{equation}
	\Diff_\rho(M) = \{\varphi\in\Diff(M)\mid \varphi_*\rho = \rho \}.
\end{equation}
The special case $\rho \equiv 1$ gives the subgroup of volume preserving diffeomorphisms denoted by $\SDiff(M)$. In general
$\varphi\in\Diff_\rho(M)$ implies that $\varphi$ is \emph{mass preserving with respect to $\rho$}.
In particular, if $\Omega \subset M$ then
\begin{equation}
	\Mass_{\rho}(\Omega) = \Mass_\rho(\varphi(\Omega)).
\end{equation}

The point of diffeomorphic density registration is to select a template density $\rho_0\in\Dens(M)$ and then generate new densities by acting on $\rho_0$ by diffeomorphisms.
In our framework we shall mostly use left action (by pushforward), but analogous results are also valid for the right action (by pullback).
One may ask `Which densities can be reached by acting on $\rho_0$ by diffeomorphisms?'
In other words find the range of the mapping
\begin{equation}
	\Diff(M)\ni \varphi \mapsto \varphi_*\rho_0 .
\end{equation}
In the language of group theory it is called the $\Diff(M)$-orbit of $\rho_0$.
This question was answered in 1965 by \citet{Mo1965} for compact manifolds: the result is that the $\Diff(M)$-orbit of $\rho_0$ consists of all densities with the same total mass as $\rho_0$. This results has been extended to non-compact manifolds \cite{greene1979diffeomorphisms} and manifolds  with boundary \cite{banyaga1974formes}. For simplicity we will only formulate the result in the compact case:



\begin{lemma}[\citet{Mo1965}]\label{lem:moser}
	Given $\rho_0,\rho_1\in \Dens(M)$, where $M$ is a compact manifold without boundary.  There exists $\varphi\in\Diff(M)$ such that $\varphi_*\rho_0 = \rho_1$ if and only if
	\begin{equation}
		\Mass_{\rho_1}(M) = \Mass_{\rho_0}(M).
	\end{equation}	
	The diffeomorphism $\varphi$ is unique up to right composition with elements in $\Diff_{\rho_0}(M)$, or, equivalently, up to left composition with elements in $\Diff_{\rho_1}(M)$.
\end{lemma}

Since the total mass of a density is a positive real number, it follows from Moser's result that the set of $\Diff(M)$-orbits in $\Dens(M)$ can be identified with $\RR_+$.
From a geometric point of view, this gives a \emph{fibration} of $\Dens(M)$ as a fiber bundle over $\RR_+$ where each fiber corresponds to a $\Diff(M)$-orbit.
In turn, Moser's result also tells us that \emph{each orbit in itself} is the base of a principal bundle fibration of $\Diff(M)$.
For example, the $\rho_0$-orbit can be identified with the quotient $\Diff(M)/\Diff_{\rho_0}(M)$ through the projection
\begin{equation}
	\pi\colon \varphi \mapsto \varphi_*\rho_0.
\end{equation}
See references~\cite{Mo2015,BaJoMo2015} for more details.


\begin{remark}
	A consequence of the simple orbit structure of $\Dens(M)$ is that one can immediately check if the registration problem can be solved exactly by comparing the total mass of $\rho_0$ and $\rho_1$.
	Furthermore, there is a natural projection from $\Dens(M)$ to any orbit simply by scaling by the total mass.

	In diffeomorphic image registration, where the action on an image is given by composition with a diffeomorphisms, the $\Diff(M)$-orbits are much more complicated. Indeed two generic images almost never belong to the same orbit. The problem of projecting from one orbit to another is ill-posed. On the other hand, because of the principal bundle structure of the space of densities, the exact registration problem of two densities with equal mass is well posed and has a complete geometric interpretation which we will exploit to develop efficient numerical algorithms.  
\end{remark}


\subsection{ $\alpha$-actions.}
The above mathematical development of diffeomorphisms acting on densities can be further generalized. 
We generalize the action of the diffeomorphism group by parametrizing the action by a positive constant $\alpha$ and define the $\alpha$-action as follows:
The group of diffeomorphisms $\Diff(\Omega)$ acts from the left on densities by the \emph{$\alpha$-action} via:
\begin{equation}
(\varphi,\rho)\mapsto \varphi_{\alpha*}\rho \coloneqq |D\varphi^{-1}|^\alpha\,\rho\circ\varphi^{-1},
\end{equation}
where $|D\varphi|$ denotes the Jacobian determinant of $\varphi$.

\begin{remark}
One of the theoretical motivation to study  $\alpha$-density action is that it allows to approximate the standard action of $\Diff(M)$ on functions given by composition: formally  $\lim_{ \alpha \rightarrow 0} \varphi_{\alpha*} I =  I\circ\varphi^{-1}$. 

From a practical point of view the motivation for the $\alpha$-action stems from the fact that CT images do not transform exactly as densities. We will  see in section~\ref{application_thoracic} that for the application of density matching of thoracic image registration, the  CT images of the lung behave as $\alpha$-densities for $\alpha<1$.
\end{remark}

Analogous to the standard mass we define the \emph{p-mass} of a subset $\Omega\subset M$ with respect to $\rho\in\Densalpha(M)$ by $$p-\Mass_\rho(\Omega) = \int_\Omega \rho^{p}\,\vol.$$
With this definition we immediately obtain the analogue of Lemma~\ref{lem:moser} for the $\alpha$-action and thus also a similar principal fiber bundle picture:
\begin{lemma}\label{lem:moser:alpha}
	Given $\rho_0,\rho_1\in \Dens(M)$, where $M$ is a compact manifold.  There exists $\varphi\in\Diff(M)$ such that $\varphi_{\alpha*}\rho_0 = \rho_1$ if and only if
	\begin{equation}
		 \frac{1}{\alpha}-\Mass_{\rho_1}(M) =  \frac{1}{\alpha}-\Mass_{\rho_0}(M).
	\end{equation}	
\end{lemma}



%
%

\section{Diffeomorphic Density registration} 
In this part we will describe a general (Riemannian) approach to diffeomorphic density registration, i.e., the problem of finding an optimal diffeomorphism $\varphi$ that transports a $\alpha$-density $\rho_0$ (source)  to a  $\alpha$-density $\rho_1$ (target). 
By Moser's result, c.f. Lemma \ref{lem:moser} and Lemma 
\ref{lem:moser:alpha} there always exists an infinite dimensional set of solutions (diffeomorphisms) to this problem. Thus the main difficulty lies in the \emph{solution selection}. Towards this aim we introduce the regularized \emph{exact $\alpha$-density registration}-problem:

\begin{tcolorbox}
Given a source density $\rho_0$ and a target density $\rho_1$  of the same total $\frac{1}{\alpha}$-mass, find a diffeomorphisms $\varphi$ that minimizes
\begin{equation}\label{eq:exact_registration}
\mathcal R(\varphi)\quad \text{under the constraint}\quad \varphi_{\alpha*}\rho_0=\rho_1\;.
\end{equation}
Here, $\mathcal R(\varphi)$ is a regularization term.
\end{tcolorbox}

\begin{remark}
Note that we have formulated the registration constraint using the left action of the diffeomorphism group, i.e., $\varphi_{\alpha*}\rho_0=\rho_1$. 
A different approach is to use the right action of $\Diff(M)$, which yields to the constraint $\varphi^{\alpha*}\rho_1=\rho_0$. These two approaches are conceptually different,
as we aim to move the source to target using the left action, while one moves the target to source using the right action. The resulting optimal deformations are however equal, if  the regularization term satisfies $\mathcal R(\varphi)=\mathcal R(\varphi^{-1})$.
\end{remark}

In the later sections we will introduce several choices for 
$\mathcal R$ and discuss their theoretical and practical properties. 
In general one aims to construct regularization terms such that the corresponding registration problem has the following desirable properties:
\begin{enumerate}
\item Theoretical results on existence and uniqueness of solutions;
\item Fast and stable numerical computations of the minimizers;
\item Meaningful optimal deformations.
\end{enumerate}
Note, that the notion of \emph{meaningful} will depend highly on the specific application.

In practice one is sometimes not interested to enforce the constraint, but is rather interested in a relaxed version of the above problem. Thus we introduce the \emph{inexact density registration}-problem: 
\begin{tcolorbox}
Given a source density $\rho_0$ and a target density $\rho_1$, find a diffeomorphisms $\varphi$ that minimizes
\begin{equation}
\mathcal E(\varphi) = \lambda\; d(\varphi_{\alpha*}\rho_0,\rho_1)+\mathcal R(\varphi)\;.
\end{equation}
Here $\lambda>0$ is a scaling parameter, $d(\cdot,\cdot)$ is a distance on the space of densities (the similarity measure) and $\mathcal R(\varphi)$ is a regularization term as before.
\end{tcolorbox}
\begin{remark}
Note, that we do not require the densities to have the same $\frac{1}{\alpha}$-mass in the inexact density matching framework. For densities that have the same $\frac{1}{\alpha}$-mass one can retrieve the exact registration problem by considering the inexact registration problem for $\lambda \to \infty$. 
\end{remark}

On the space of probability densities there exists a canonical Riemannian metric, the Fisher-Rao metric, which allows for explicit formulas of the corresponding geodesic distance: it is given by the (spherical) Hellinger distance. For the purpose of this book chapter we will often use this distance functional as a similarity measure. We will however discuss several choices of different regularization terms, which will be the topic of the next sections.

\section{Density registration in the LDDMM-framework}\label{LDDMM_density}
The LDDMM-framework is based on the idea of using a \emph{right-invariant} metric on the diffeomorphism group to define the regularity measure, i.e.,
\begin{equation}
\mathcal R(\varphi)=\operatorname{dist}(\operatorname{id},\varphi)\;,
\end{equation}
where $\operatorname{dist}(\cdot,\cdot)$ denotes the geodesic distance of a right invariant metric on $\Diff(M)$. 

\begin{remark}
In the standard presentation of the LDDMM framework right-invariant metrics on $\Diff(M)$ are usually defined using the theory of Reproducing Kernel Hilbert Spaces (RKHS).
We will follow a slightly different approach and equip the whole group of diffeomorphisms with a weak right-invariant metric, see \cite{BrVi2014} for a comparison of these two approaches.   
\end{remark}

From here on we assume that $M$ is equipped with a smooth Riemannian metric $g$ with volume density $\mu$.
To define a right invariant metric on  $\Diff(M)$ we introduce the so-called \emph{inertia operator} 
$A\colon \mathfrak X(M)\to \mathfrak X(M)$, 
where $\mathfrak X(M)$, the set of smooth vector fields, is the Lie-Algebra of $\Diff(M)$.
We will assume that $A$ is a strictly positive, elliptic, differential operator, that is self adjoint with respect to the $L^2$ inner product on $\mathfrak X(M)$.
For the sake of simplicity we will only consider operators $A$, that are defined via powers of the Laplacian of the Riemannian metric $g$, i.e. 
we will only consider operators of the form
\begin{equation}\label{inertia_operator}
A=(1-\Delta_g)^k
\end{equation}
for some integer $k$.
Here $\Delta_g$ denotes the Hodge-Laplacian of the metric $g$.
Most of the results discussed below are valid for a much larger class of (pseudo) differential operators, see \cite{BaJoMo2017}. Any such $A$ defines an inner product $G_{\operatorname{\id}}$ on  $\mathfrak X(M)$ via
\begin{equation}
 G_{\operatorname{\id}}(X,Y)=\int_M g\left(A X,  Y\right)\, \mu\,.
\end{equation}
where $\mu$ denotes the induced volume density of $g$.
We can extend this to a right-invariant metric on $\Diff(M)$ by right-translation:
\begin{equation}\label{eq:Gmet_DiffM}
G_{\varphi}(h,k)= G_{\operatorname{\id}}(h\circ\varphi^{-1},k\circ\varphi^{-1})= \int_M g\left(A (h\circ\varphi^{-1}),  k\circ\varphi^{-1}\right)\, \mu\;.
\end{equation}
For an overwiew on right invariant metrics on diffeomorphism groups we refer to the articles \cite{MP2010,BrVi2014,bauer2016uniqueness,BEK2015}. 

In this framework the exact density registration problem reads as:
\begin{tcolorbox}
Given a source density $\rho_0$ and a target density $\rho_1$ find a diffeomorphisms $\varphi$ that minimizes
\begin{equation}
\operatorname{dist}(\operatorname{id},\varphi)\qquad \text{such that } \varphi_{\alpha*}\rho_0=\rho_1\;.
\end{equation}
where  $\operatorname{dist}(\operatorname{id},\varphi)$ is the geodesic distance on $\Diff(M)$ of the metric \eqref{eq:Gmet_DiffM}.
\end{tcolorbox}
Using this particular regularization term provides an intuitive interpretation of the  solution selection: one aims to find the transformation that is as close as possible to the identity under the constraint that it transports the source density to the target density.
 
In the following Theorem we present a summary of the geometric picture, that underlies the exact registration problem. 
To keep the presentation simple we will only consider the case $\alpha=1$, i.e., the standard density action. A similar result can be obtained for general $\alpha$.

Let $\pi$ be the projection 
\begin{equation}
\pi: \Diff(M)\rightarrow \Dens(M)\simeq \Diff_{\rho}(M)\backslash\Diff(M)
\end{equation}
that is induced  by the left action of the diffeomorphism group, c.f. Lemma~\ref{lem:moser}. Then we have:
\begin{theorem}
Let $G$ be a right-invariant metric on $\Diff(M)$ of the form~\eqref{eq:Gmet_DiffM} with inertia operator $A$ as in~\eqref{inertia_operator}.
Then there exists a unique metric $\bar G_{\rho}$ on $\Dens(M)$ such that the projection $\pi$ is a Riemannian submersion.
The order of the induced metric $\bar G$ on $\Dens(M)$ is  $k-1$, where $k$ is the order of the metric $G$. 
\end{theorem}
A direct consequence of the Riemannian submersion picture is the following characterization of the solutions of the exact density registration problem:
\begin{corollary}
Let $\rho(t)$, for $t\in[0,1]$, be a minimizing geodesic connecting the given densities $\rho_0$ (source) and $\rho_1$ (target). Then the solution of the exact registration problem is  
given by the endpoint $\varphi(1)$ of the horizontal lift of the geodesic $\rho(t)$. 
\end{corollary}
\begin{remark}
The above result describes an intriguing geometric interpretation of the solutions of the exact registration problem. 
Its applicability is however limited to cases where there exist an explicit solutions for the geodesic boundary value problem on the space of probability densities with respect to the metric $\bar G$. 
To our knowledge the only such example is the so-called Optimal Information Transport setting, which we will discuss in the next section. In the general case the solution of the exact
density registration problem requires one to solve the horizontal  geodesic boundary value problem on the group of diffeomorphisms, which is connected to the 
solution of a nonlinear PDE, the EPDiff equation. Various algorithms have been proposed for numerically solving the optimization problems~\cite{Beg2005,Younes2009S40,FX2012}.
\end{remark}


\section{Optimal Information Transport}
In this section we describe an explicit way of solving the exact density registration problem. 
The framework  in this section has been previously developed for random sampling from non-uniform arbitrary distributions  \cite{BaJoMo2017b}. For simplicity we will restrict ourself to the standard density action, i.e., $\alpha=1$. However, all the algorithms are easily generalized to general $\alpha$.
The specific setting uses deep geometric connections between the Fisher--Rao metric on the space of probability densities and a special right-invariant metric on the group of diffeomorphisms. 

\begin{definition}
	The \emph{Fisher--Rao metric} is the Riemannian metric on $\Dens(M)$ given by
	\begin{equation}
		G_\rho(\dot\rho,\dot\rho) = \frac{1}{4}\int_M \frac{\dot\rho^2}{\rho}\,\vol .
	\end{equation}
\end{definition}
The main advantage of the Fisher-Rao metric is the existence of explicit formulas for the solution to the geodesic boundary value problem and thus also for the induced geodesic distance:
\begin{proposition}[\citet{Fr1991}]
	Given $\rho_0,\rho_1\in\Dens(M)$ with the same total mass, the Riemannian distance with respect to the Fisher--Rao metric is given by
	\begin{equation}\label{eq:fr_distance}
		d_F(\rho_0,\rho_1) = \arccos\left( \int_M \sqrt{\frac{\rho_1}{\rho_0}}\rho_0 \right)
	\end{equation}
	Furthermore, the geodesic between $\rho_0$ and $\rho_1$ is given by
	\begin{equation}\label{eq:fr_geodesics}
		\rho(t) = \left(\frac{\sin((1-t)\theta)}{\sin\theta}+ \frac{\sin(t\theta)}{\sin\theta}\sqrt{\frac{\rho_1}{\rho_0}}\right)^2\rho_0
	\end{equation}
	where $\theta = d_F(\rho_0,\rho_1)$.
\end{proposition}

Using the formula \eqref{eq:fr_geodesics} for geodesics we shall construct an almost explicit algorithm for solving an exact density registration problem of the form in \eqref{eq:exact_registration}.
To this end we need to introduce a suitable regularization term. 
As in the LDDMM framework (see  section~\ref{LDDMM_density}) we shall chose it as distance to the identity with respect to a right-invariant Riemannian metric on $\Diff(M)$.
However, in order to exploit the explicit formula~\eqref{eq:fr_geodesics} the right-invariant metric needs to communicate with the Fisher--Rao metric, as we now explain.

\begin{definition}
	The \emph{information metric} is the right-invariant Riemannian metric on $\Diff(M)$ given (at the identity) by
	\begin{equation}
		\bar G_\id(u,v) = -\int_M \pair{\Delta u,v}\,\vol + 
		 \sum_{i=1}^k \int_M \pair{ u,\xi_i}\,\vol  \int_M \pair{v,\xi_i}\,\vol
	\end{equation}
	where $\Delta u$ denotes the Laplace-de Rham operator lifted to vector fields and where $\xi_1,\ldots,\xi_k$ is a basis of the harmonic fields on $M$. The Riemannian distance corresponding to $\bar G$ is denoted $d_I(\cdot,\cdot)$. Because of the Hodge decomposition theorem, the metric is independent of the choice of orthonormal basis for the harmonic fields.

Building on work by Khesin, Lenells, Misiolek, and
Preston \cite{KhLeMiPr2013}, Modin \cite{Mo2015} showed that the metric $\bar G $ descends to the Fisher-Rao metric on the space of densities. This fundamental property will serve as the basis for our algorithms. \end{definition}

We are now ready to formulate our special density registration problem, called the \emph{optimal information transport problem}:

\begin{tcolorbox}
	\textbf{Optimal information transport (OIT)}
	\smallskip

	Given $\rho_0,\rho_1\in\Dens(M)$, and a Riemannian metric on $M$ with volume form $\rho_0$, find a diffeomorphism $\varphi$ that minimize
	\begin{equation}\label{eq:oit_registration}
		\mathcal E(\varphi) = d_I(\id,\varphi) = d_F(\rho_0,\varphi_*\rho_0)
	\end{equation}
	under the constraint $\varphi_*\rho_0=\rho_1$.
\end{tcolorbox}

In general, the formula for $d_I(\id,\varphi)$ is not available explicitly; one would have to solve a nonlinear PDE (the EPDiff equation).
However, because of the special relation between $d_I$ and $d_F$ we have the following result, which is the key to an efficient algorithm.

\begin{theorem}[\cite{Mo2015,BaJoMo2015}]\label{thm:oit}
The OIT problem has a unique solution. 
That is, there is a unique diffeomorphism $\varphi\in \Diff(M)$ minimizing $d_I(\id,\varphi)$ under the constraint $\varphi_*\rho_0=\rho_1$.
The solution is explicitly given by $\varphi(1)$, where $\varphi(t)$ is the solution to the problem
\begin{equation}\label{eq:veq}
\begin{split}
	\Delta f(t) &= \frac{\dot\rho(t)}{\rho(t)}\circ \varphi(t), \\
	v(t) &= \nabla(f(t)), \\
	\frac{d}{d t}\varphi(t)^{-1} &= v(t)\circ\varphi(t)^{-1}, \quad \varphi(0) = \id 
\end{split}
\end{equation}
and $\rho(t)$ is the Fisher--Rao geodesic connecting $\rho_0$ and $\rho_1$
\begin{equation}\label{eq:fisher_rao_geodesics}
	\rho(t)= \left( 
			\frac{\sin\left((1-t)\theta\right)}{\sin\theta} + \frac{\sin\left(t\theta\right)}{\sin\theta}\sqrt{\frac{\rho_1}{\rho_0}}
		\right)^2 \rho_0,
		\quad \cos\theta = \int_M   \sqrt{\frac{\rho_1}{\rho_0}}\; \rho_0  \,.
\end{equation}
\end{theorem}

Based on Theorem~\ref{thm:oit} we now give a semi-explicit algorithm for numerical computation of the solution to the optimal information transport problem.
The algorithm assumes that we have a numerical way to represent functions, vector fields, and diffeomorphisms on~$M$, and numerical methods for
\begin{itemize}
\item composing functions and vector fields with diffeomorphisms,
\item computing the nablaient of functions, and
\item computing solutions to Poisson's equation on $M$.
\end{itemize}

\begin{tcolorbox}
	\textbf{Numerical algorithm for optimal information transport}
	\begin{enumerate}
	\item Choose a step size $\varepsilon = 1/K$ for some positive integer~$K$ and calculate the Fisher-Rao geodesic $\rho(t)$ and its derivative $\dot \rho(t)$ at all time points $t_k=\frac{k}{K}$ using equation~\eqref{eq:fisher_rao_geodesics}.
		\item 
		Initialize $\varphi_0 = \id$. 
		Set $k\leftarrow 0$.

		\item Compute $s_k = \frac{\dot \rho(t_k)}{\rho(t_k)}\circ\varphi_{k}$ and solve the Poisson equation
		\begin{equation}
			\Delta f_k= s_k.
		\end{equation}

		\item Compute the gradient vector field $v_k = \nabla f_k$.

		\item Construct approximations $\psi_k$ to $\exp(-\varepsilon v_k)$, for example 
		\begin{equation}
			\psi_k = \id - \varepsilon v_k.		
		\end{equation}

		\item Update the diffeomorphism\footnote{If needed, one may also compute the inverse by $\varphi_{k+1}^{-1} = \varphi_k^{-1} + \varepsilon v\circ\varphi_k^{-1}$.}
		\begin{equation}
			\varphi_{k+1} = \varphi_k\circ\psi_k .
		\end{equation}
		
		\item Set $k \leftarrow k+1$ and continue from step 3 unless $k=K$.
	\end{enumerate}
\end{tcolorbox}

Although it is possible to use optimal information transport and the algorithm above for medical image registration problems, the results so obtained are typically not satisfactory; the diffeomorphism obtained tends to compress and expand matter instead of moving it (see example in \cite[Sec.~4.2]{BaJoMo2015}).
Another problem is that the source and target densities are required to be strictly positive, which is typically not the case for medical images.
However, in application where either the source or the target density is uniform (with respect to the natural Riemannian structure of the manifold at hand), the OIT approach can be very competitive.

\subsection{Application: random sampling from non-uniform distribution}
In this section we describe an application of OIT to random sampling from non-uniform distributions, i.e., the following problem.

\begin{tcolorbox}
\textbf{Random sampling problem}
\smallskip

Let $\rho_1\in\Dens(M)$.
Generate $N$ random samples from the probability distribution $\rho_1$. 
\end{tcolorbox}

The classic approach to sample from a probability distribution on a higher dimensional space is to use Markov Chain Monte Carlo (MCMC) methods, for example the Metropolis--Hastings algorithm~\cite{Ha1970}. 
An alternative idea is to use diffeomorphic density registration between the density $\rho_1$ and a standard density $\rho_0$ from which samples can be drawn easily.
Indeed, one can then draw samples from $\rho_0$ and transform them via the computed diffeomorphism to generate samples from $\rho_1$.
A benefit of transport-based methods over traditional MCMC methods is cheap computation of additional samples;
it amounts to drawing uniform samples and then evaluating the transformation.
On the other hand, transport-based methods scale poorly with increasing dimensionality of $M$, contrary to MCMC.

Moselhy and Marzouk \cite{MoMa2012} and Reich \cite{Re2013} proposed to use optimal mass transport (OMT) to construct the desired diffeomorphism $\varphi$, thereby enforcing $\varphi = \nabla c$ for some convex function~$c$. 
The OMT approach implies solving, in one form or another, the heavily non-linear Monge--Ampere equation for $c$.
A survey of the OMT approach to random sampling is given by Marzouk \emph{et.~al.}~\cite{MaMoPaSp2016}.
Using OIT instead of OMT, the problem simplifies significantly, as the OIT-algorithm above only involves solving linear Poisson problems.

\begin{figure}[http]
	\includegraphics[height=0.3\textheight]{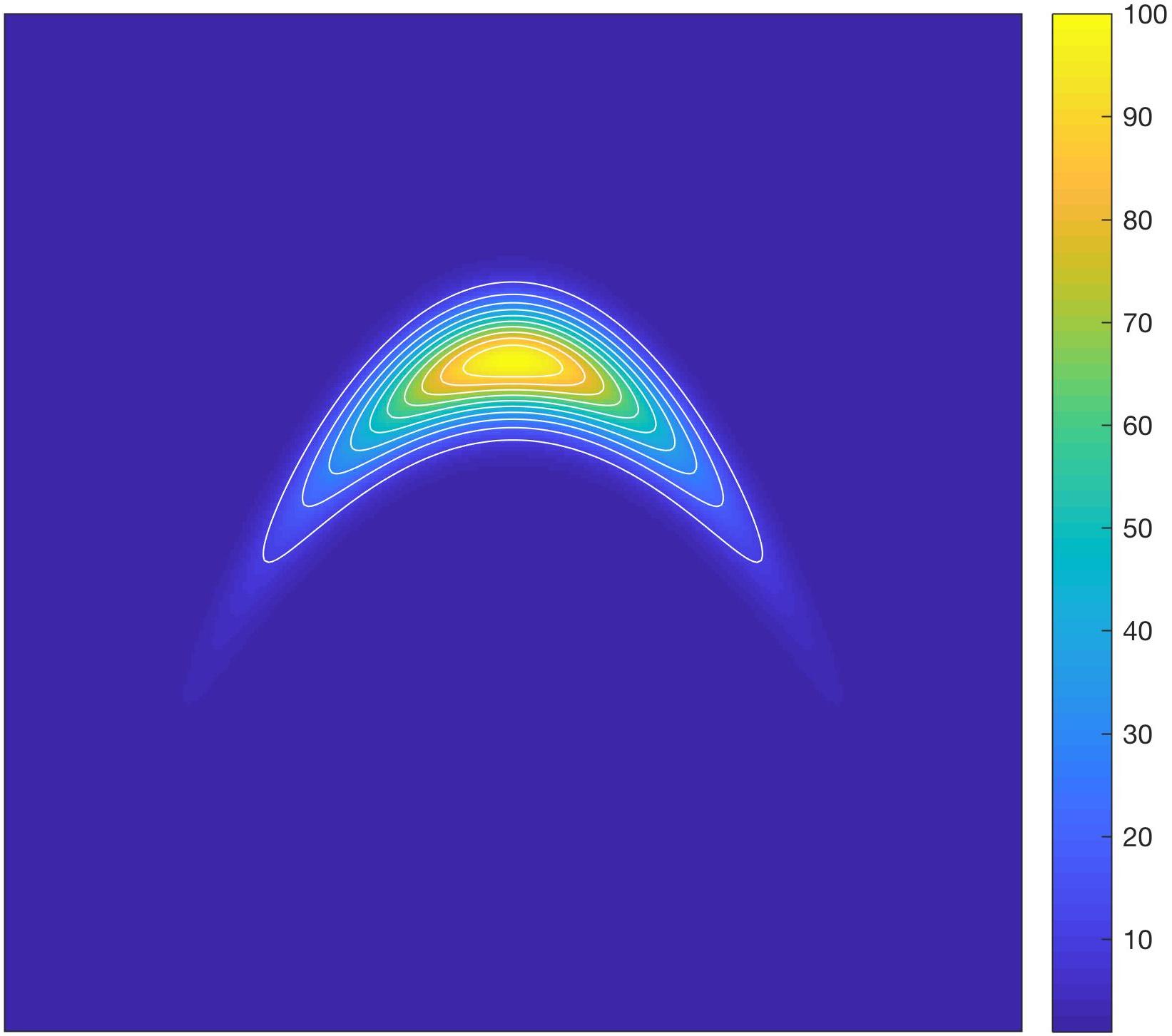}
	\includegraphics[height=0.3\textheight]{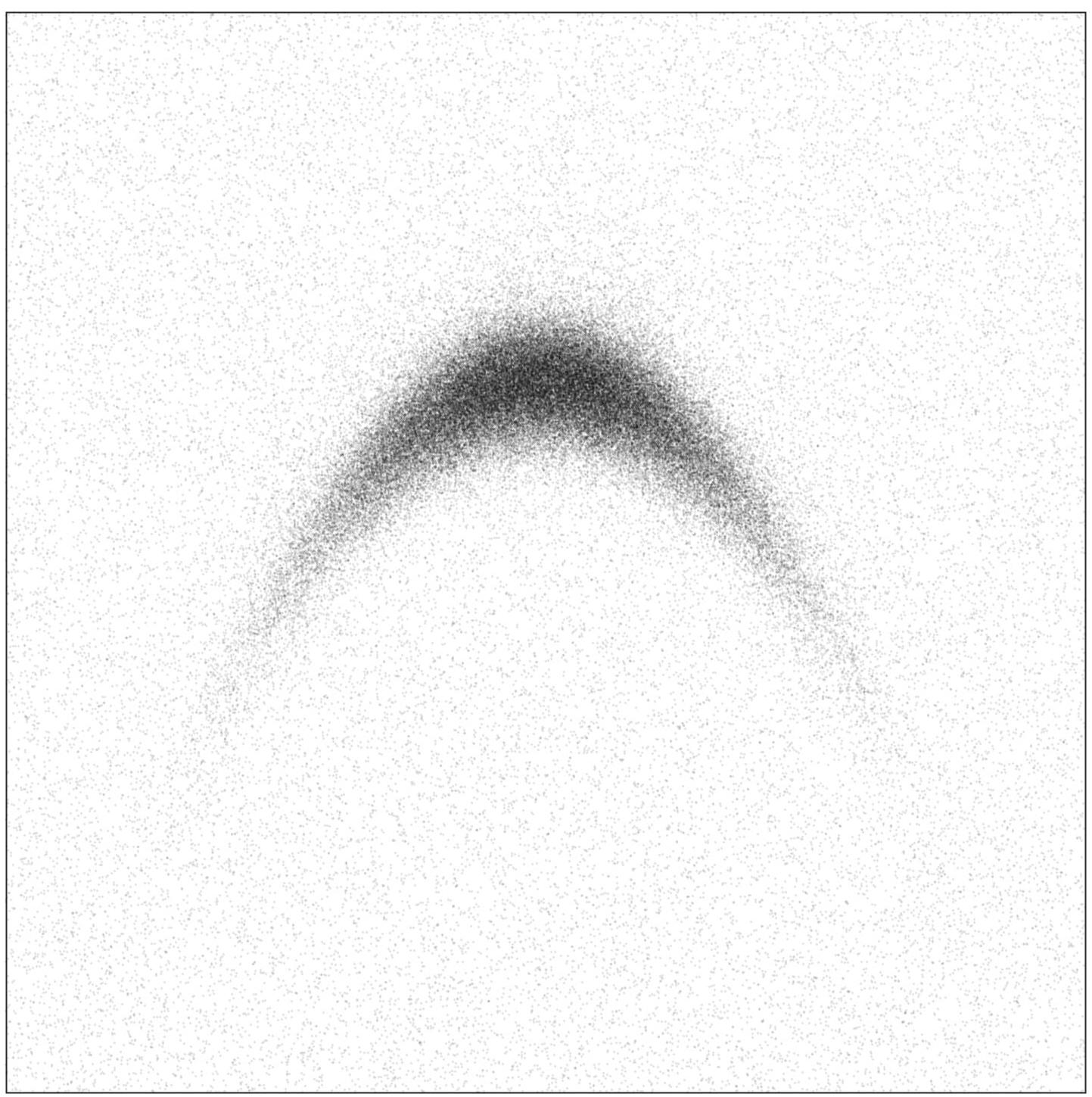}
	\caption{Application of OIT to random sampling: (left) The probability density $\rho$. The maximal density ratio is 100.
	(right) $10^5$ samples from $\rho$ calculated using our OIT based random sampling algorithm.}
	\label{fig:example1_samples}
\end{figure}

As a specific example, consider consider $M=\mathbb{T}^2 \simeq (\mathbb{R}/2\pi\mathbb{Z})^2$ with distribution defined in Cartesian coordinates $x,y\in [-\pi,\pi)$ by
\begin{equation}\label{distribution}
	\rho \sim 3\exp(-x^2 - 10(y-x^2/2+1)^2) +1/10,
\end{equation}
normalized so that the ratio between the maximum and mimimum of $\rho$ is 100.
The resulting density is depicted in Fig.~\ref{fig:example1_samples}~(left).

We draw $10^5$ samples from this distribution using a MATLAB implementation of our algorithm, available under MIT license at
\begin{center}
	\texttt{\url{https://github.com/kmodin/oit-random}}
\end{center}

\begin{figure}[http]
	\begin{center}
		\includegraphics[height=0.48\textheight]{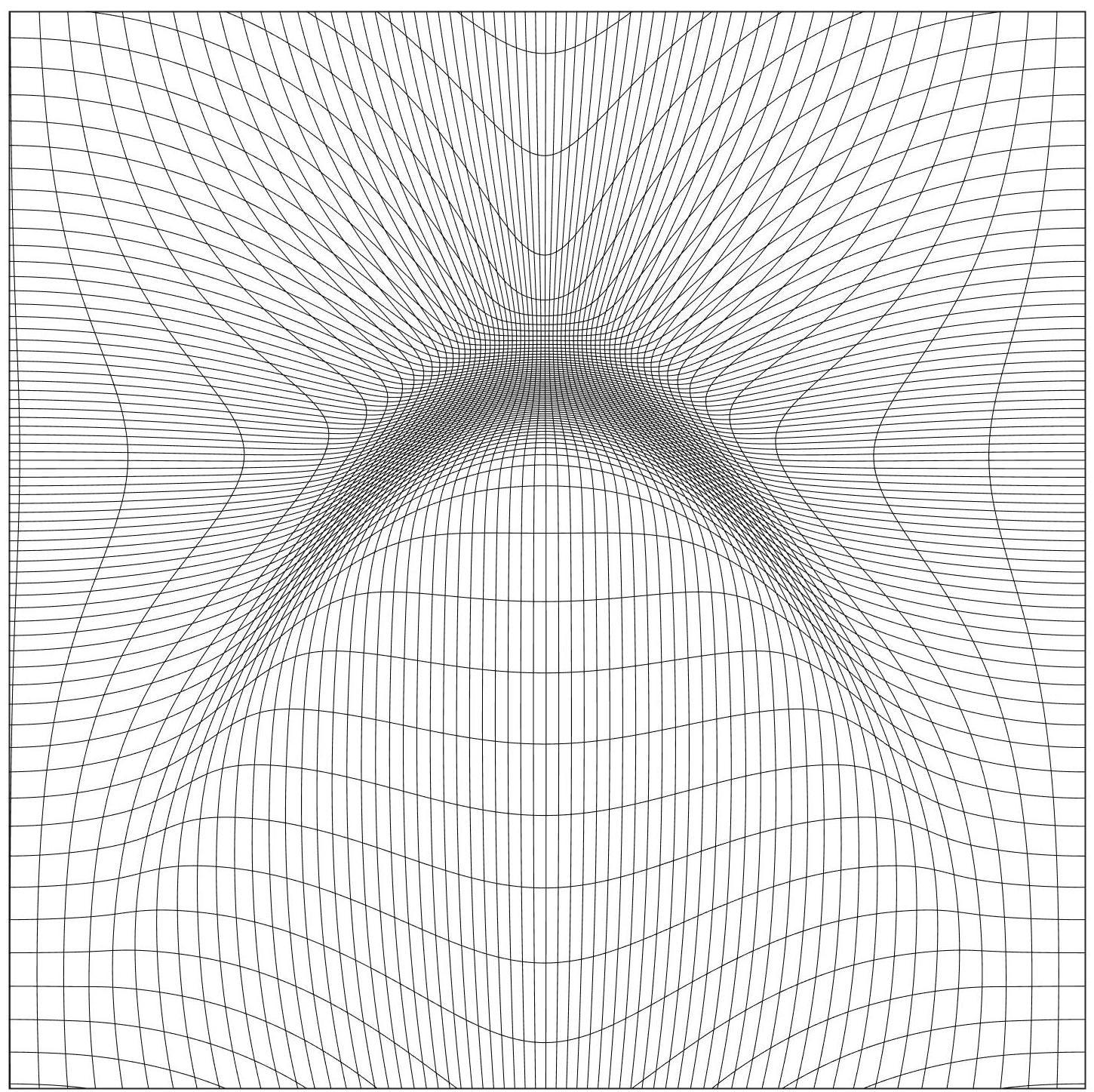}
	\end{center}
	\caption{Application of OIT to random sampling:
	The computed diffeomorphism $\varphi_K$ shown as a warp of the uniform $256\times 256$ mesh (every 4th mesh-line is shown). 
	Notice that the warp is periodic.
	The ratio between the largest and smallest warped volumes is 100.
	}
	\label{fig:example1_phiinv}
\end{figure}

The implementation can be summarized as follows.
To solve the Poisson problem we discretize the torus by a $256\times 256$ mesh and use the fast Fourier transform (FFT) to invert the Laplacian. 
We use 100 time steps.
The resulting diffeomorphism is shown as a mesh warp in Fig.~\ref{fig:example1_phiinv}.
We then draw $10^5$ uniform samples on $[-\pi,\pi]^2$ and apply the diffeomorphism on each sample (applying the diffeomorphism corresponds to interpolation on the warped mesh).
The resulting random samples are depicted in Fig.~\ref{fig:example1_samples}~(right).
To draw new samples is very efficient.
For example, another $10^7$ samples can be drawn in less than a second.

\section{A gradient flow approach}
In the optimal information transport described in the previous section the fundamental restriction is that the volume form of Riemannian metric of the base manifold is compatible with the density being transformed in that it has to be conformally related to the source density $\rho_0$. In most medical imaging application this modeling assumption is not applicable. In this section we will develop more general algorithms that relax the requirement for the metric to be compatible with the densities to be registered. In this section we will consider the natural extension of the Fisher-Rao metric to the space of all densities and the case when $\vol(\Omega)=\infty$, for which it is given by
\begin{equation}\label{eq:FR_distance_infinite}
	d_F^{2}(I_0\,dx,I_1\, dx)= \int_{\Omega} (\sqrt{I_0}-\sqrt{I_1})^{2} dx \, .
\end{equation}
Notice that $d_F^{2}(\cdot,\cdot)$ in this case is the \emph{Hellinger distance}.
For details, see~\cite{BaJoMo2015}.

The Fisher--Rao metric is the unique Riemannian metric on the space of probability densities that is invariant under the action of the diffeomorphism group~\cite{bauer2016uniqueness,AJLS2014}.
This invariance property extends to the induced distance function, so
\begin{equation}\label{eq:invariance}
	d_F^{2}(I_0\,dx,I_1\,dx)= d_F^{2}(\varphi_*(I_0\,dx),\varphi_*(I_1\,dx)) \qquad \forall \varphi \in \Diff(\Omega)\;.
\end{equation}

\begin{figure}[htp]
	\centering
	\includegraphics[width=0.7\textwidth]{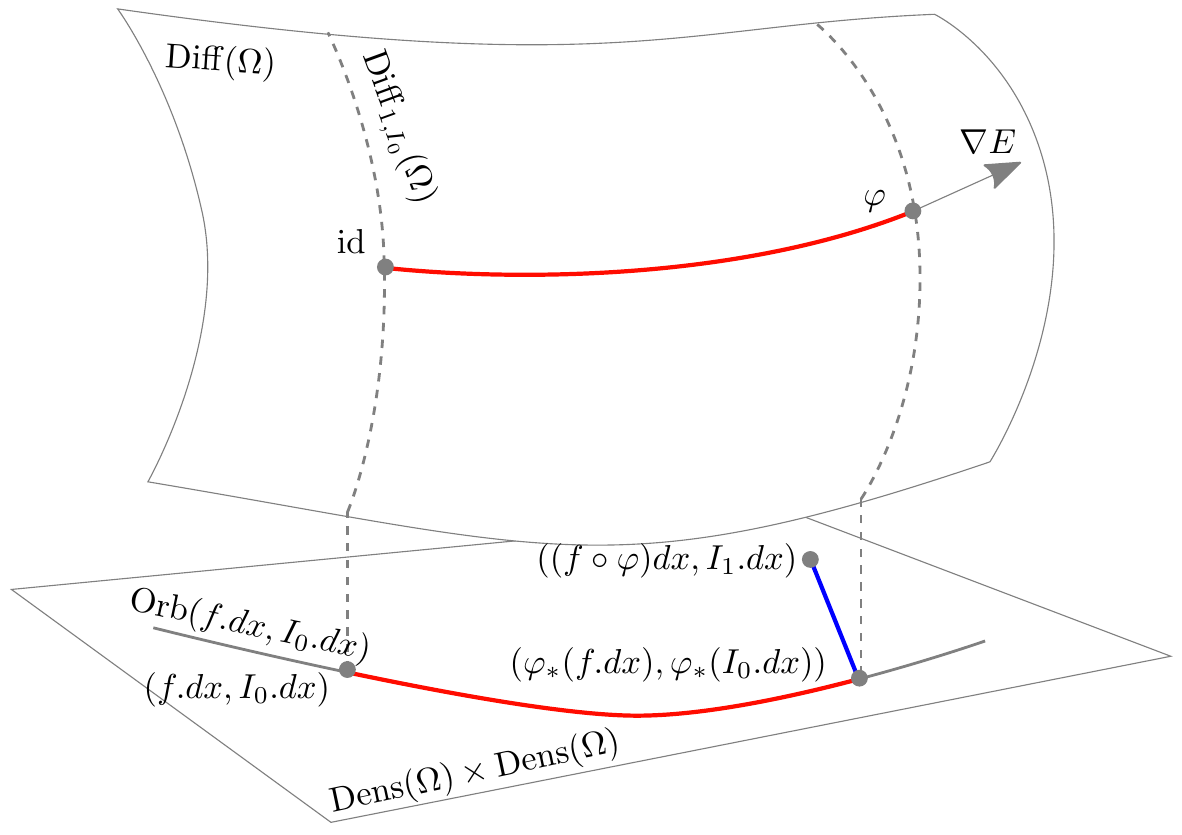}
	\caption{Illustration of the geometry associated with the density registration problem.
	The gradient flow on $\operatorname{Diff}(\Omega)$ descends to a gradient flow on the orbit $\mathrm{Orb}(f\, dx,I_0 \, dx)$.
	While constrained to $\mathrm{Orb}(f\,dx,I_0\,dx)\subset \operatorname{Dens}(\Omega)\times\operatorname{Dens}(\Omega)$, this flow strives to minimize the product Fisher-Rao distance to $((f\circ\varphi)\,dx,I_1\,dx)$.
	}
	\label{Fig1}
\end{figure}

Motivated by the aforementioned properties, we develop a weighted diffeomorphic registration algorithm for registration two density images.
The algorithm is based on the Sobolev $H^{1}$ gradient flow on the space of diffeomorphisms that minimizes the energy functional
\begin{equation}\label{eq:E}
E(\varphi) = d^{2}_F( \varphi_*(f\, dx), (f \circ \varphi^{-1}) dx) + d^{2}_F( \varphi_*(I_0\, dx), I_1\, dx)) .
\end{equation}
This energy functional is only a slight modification of the energy functional studied in~\cite{BaJoMo2015}.
Indeed, if $f$ in the above equation is a constant $\sigma>0$, then~\eqref{eq:E} reduces to the energy functional of Bauer, Joshi, and Modin~\cite[\S\!~5.1]{BaJoMo2015}.
Moreover, the geometry described in~\cite[\S\!~5.3]{BaJoMo2015} is valid also for the functional~\eqref{eq:E}, and, consequently, the algorithm developed in~\cite[\S\!~5.2]{BaJoMo2015} can be used also for minimizing~\eqref{eq:E}.
There the authors view the energy functional as a constrained minimization problem on the product space $\Dens(\Omega)\times\Dens(\Omega)$ equipped with the
product distance, cf. Fig~\ref{Fig1} and ~\cite[\S\!~5]{BaJoMo2015} for details on the resulting geometric picture. Related work on  diffeomorphic  density registration using the Fisher Rao metric can be found in \cite{STFV2013,SHV2013}.

Using the invariance property of the Fisher-Rao metric and assuming infinite volume, the main optimization problem associated with the energy functional~\eqref{eq:E} is the following.

	\noindent Given densities $I_0\,dx$, $I_1\,dx$, and $f\, dx$, find $\varphi\in\operatorname{Diff}(\Omega)$ minimizing
	\begin{equation}\label{eq:main_problem}
		E(\varphi) = \underbrace{\int_{\Omega} (\sqrt{|D\varphi^{-1}|}-1)^2\,f\circ\varphi^{-1} \,dx}_{E_1(\varphi)} +
		\underbrace{\int_{\Omega}\Big(\sqrt{|D\varphi^{-1}|I_0\circ\varphi^{-1}}- \sqrt{I_1} \Big)^2\,dx}_{E_2(\varphi)} \;.
	\end{equation}

The invariance of the Fisher-Rao distance can be seen with a simple change of variables $x \mapsto \varphi(y)$, $dx \mapsto |D\varphi| dy$, and $|D\varphi^{-1}| \mapsto \frac{1}{|D\varphi|}$.
Then, Equation \ref{eq:main_problem} becomes
\begin{align}
  \label{eq:main_problemCV}
 E(\varphi) = \int_{\Omega} (1-\sqrt{|D\varphi|})^2\,f \,dy +
 \int_{\Omega}\Big(\sqrt{I_0}- \sqrt{|D\varphi|I_1\circ\varphi} \Big)^2\,dy \;.
\end{align}
To better understand the energy functional $E(\varphi)$ we consider the two terms separately.
The first term $E_1(\varphi)$ is a \emph{regularity measure} for the transformation.
It penalizes the deviation of the diffeomorphism $\varphi$ from being volume preserving.
The density $f\,dx$ acts as a weighting on the domain $\Omega$.
That is, change of volume (compression and expansion of the transformation $\varphi$) is penalized more in regions of $\Omega$ where $f$ is large.
The second term $E_2(\varphi)$  penalizes \emph{dissimilarity} between $I_0\, dx$ and $\varphi^*(I_1\, dx)$.
It is the Fisher--Rao distance between the initial density $I_0\, dx$ and the transformed target density $\varphi^*(I_1\, dx)$.
Because of the invariance~\eqref{eq:invariance} of the Fisher--Rao metric, this is the same as the Fisher--Rao distance between $I_1\, dx$ and $\varphi_*(I_0\, dx)$.

Solutions to problem~\eqref{eq:main_problem} are \emph{not} unique.
To see this, let $\Diff_{I}(\Omega)$ denote the space of all diffeomorphisms preserving the volume form $I\,dx$:
\begin{equation}
	\Diff_{I}(\Omega) = \{ \varphi\in\Diff(\Omega)\mid \; |D\varphi|\,(I\circ\varphi)=I \}	.
\end{equation}
If $\varphi$ is a minimizer of $E(\cdot)$, then $\psi\circ\varphi$ for any
\begin{equation}
	\psi\in\Diff_{1,I_0}(\Omega) \coloneqq  \Diff_{1}(\Omega)\cap \Diff_{I_0}(\Omega)
\end{equation}
is also a minimizer.
Notice that this space is not trivial.
For example, any diffeomorphism generated by a \emph{Nambu--Poisson vector field} (see~\cite{Na1999}), with $I_0$ as one of its Hamiltonians, will belong to it.
A strategy to handle the degeneracy was developed in~\cite[\S\!~5]{BaJoMo2015}: the fact that the metric is descending with respect to the $H^1$ metric on $\Diff(\Omega)$ can be used to ensure that the gradient flow is \emph{infinitesimally optimal}, i.e., always orthogonal to the null-space.
We employ the same strategy in this paper.
The corresponding geometric picture can be seen in Fig.~\ref{Fig1}.

To derive an gradient algorithm to optimize the energy functional the natural metric on the space of diffeomorphisms to use if the $H^{1}$-metric due to it's intimate link with the Fisher-Rao metric as described previously. 
The $H^{1}$-metric on the space of diffeomorphisms is defined using the Hodge laplacian on vector fields and is given by:
\begin{equation}\label{eq:H1dot}
	\GI_{\varphi}(U,V) = \int_{\Omega}\langle -\Delta u,v\rangle dx \ .
\end{equation}
Due to its connections to information geometry we also refer to this metric as \emph{information metric}.
 Let $\nabla^{\GI} E$ denote the gradient with respect to the information metric defined above. Our approach to minimize the functional of~\eqref{eq:main_problemCV} is to use a simple Euler integration of the discretization of the gradient flow:
 \begin{equation}
 {\dot \varphi } = - \nabla^{\GI} E(\varphi)
 \end{equation}
 The resulting final algorithm is order of magnitudes faster than LDDMM, since we are not required to time integrate the geodesic equations, as necessary in LDDMM~\cite{Younes2009S40}.

 In the following theorem we calculate the gradient of the energy functional:
 \begin{theorem}
  The $\GI$--gradient of the registration functional~\eqref{eq:main_problemCV} is given by
  \begin{align}
  \nabla^{\GI} E &= -\Delta^{-1} \Big( -\grad\big(f\circ \varphi^{-1} (1-\sqrt{|D\varphi^{-1}|})\big) - \nonumber \\
  & \hspace{30pt} \sqrt{|D\varphi^{-1}|\,I_0\circ\varphi^{-1}} \grad\big(\sqrt{I_1}\big) + \grad\big(\sqrt{|D\varphi^{-1}|\,I_0\circ\varphi^{-1}}\big)\sqrt{I_1} \Big)\;.
\end{align}
\label{Thm-H1Grad}
 \end{theorem}
 \begin{remark}
Notice that in the formula for $\nabla^{\GI} E$ we never need to compute $\varphi$, so in practice we only compute $\varphi^{-1}$.
We update this directly via $\varphi^{-1}(x) \mapsto \varphi^{-1}(x + \epsilon \nabla^{\GI} E)$ for some step size $\epsilon$.
\end{remark}
\begin{proof}
We first calculate the variation of the energy functional. Therefore let $\varphi_s$ be a family of diffeomorphisms parameterized by the real variable $s$, such that
\begin{equation}
  \varphi_0 = \varphi \hspace{10pt} \mathrm{and} \hspace{10pt} \dd{}{s} \Big|_{s=0} \varphi_s = v \circ \varphi.
\end{equation}
We use the following identity, as derived in \cite{hinkle2013idiff}:
\begin{align}
  \dd{}{s} \Big|_{s=0} \sqrt{|D\varphi_s|} =& \frac{1}{2}\sqrt{|D\varphi|}\mathrm{div}(v) \circ \varphi .
\end{align}
The variation of the first term of the energy functional is
\begin{align}
    \dd{}{s} \Big|_{s=0} E_1(\varphi) =& \int_\Omega f(y) (\sqrt{|D\varphi(y)|}-1)\sqrt{|D\varphi(y)|} \mathrm{div}(v) \circ \varphi(y) dy
\end{align}
We do a change of variable $y \mapsto \varphi^{-1}(x)$, $dy \mapsto |D\varphi^{-1}(x)| dx$, using the fact that $|D\varphi(y)| = \frac{1}{|D\varphi^{-1}(x)|}$;
\begin{align}
  =& \int_\Omega f\circ \varphi^{-1}(x) (1-\sqrt{|D\varphi^{-1}(x)|}) \mathrm{div}(v(x)) dx \\
  =& \left \langle f\circ \varphi^{-1} (1-\sqrt{|D\varphi^{-1}|}), \mathrm{div}(v) \right \rangle_{L^2(\RR^3)} \\
  =& -\left \langle \grad\left(f\circ \varphi^{-1} (1-\sqrt{|D\varphi^{-1}|})\right), v \right \rangle_{L^2(\RR^3)}
\end{align}
using the fact that the adjoint of the divergence is the negative gradient.
For the second term of the energy functional, we expand the square
\begin{align}
  E_2(\varphi) &= \int_\Omega  I_0(y) - 2\sqrt{I_0(y) I_1 \circ \varphi(y) |D\varphi(y)|} + I_1 \circ \varphi(y)|D\varphi(y)| dy
\end{align}
Now $\int_\Omega I_1 \circ \varphi(y)|D\varphi(y)| dy$ is constant (conservation of mass), so we only need to minimize over the middle term.
The derivative is then
\begin{align}
  \dd{}{s} \Big|_{s=0} E_2(\varphi) &= - \int_\Omega 2 \sqrt{I_0(y)} \big(\grad \sqrt{I_1}^T v\big)\circ \varphi(y) \sqrt{|D\varphi(y)|} \nonumber \\
 & \hspace{20pt}  - \sqrt{I_0(y)I_1 \circ \varphi(y) |D\varphi(y)|} \mathrm{div}(v) \circ \varphi(y) dy.
\end{align}
We do the same change of variables as before:
\begin{align}
   &= - \int_\Omega \sqrt{I_0\circ\varphi^{-1}(x)} \frac{|D\varphi^{-1}(x)|}{\sqrt{|D\varphi^{-1}(x)}|} \big( 2\grad \sqrt{I_1(x)}^T v(x) +\sqrt{I_1(x)}\mathrm{div}(v)(x) \big) \\
   &= -\left \langle 2 \sqrt{|D\varphi^{-1}|\,I_0\circ\varphi^{-1}} \grad \sqrt{I_1}, v \right \rangle_{L^2(\RR^3)} \nonumber \\
   &\hspace{50pt} - \left \langle  \sqrt{|D\varphi^{-1}|\,I_0\circ\varphi^{-1} I_1}, \mathrm{div}(v) \right \rangle_{L^2(\RR^3)} \\
   &= \left \langle - \sqrt{|D\varphi^{-1}|\,I_0\circ\varphi^{-1}} \grad \sqrt{I_1}, v \right \rangle_{L^2(\RR^3)} \nonumber \\
    &\hspace{50pt} + \left \langle \grad\left(\sqrt{|D\varphi^{-1}|\,I_0\circ\varphi^{-1}}\right)\sqrt{I_1}, v \right \rangle_{L^2(\RR^3)}.
\end{align}
From the above equations we conclude that:
\begin{align}
  \label{eq:finalgrad}
  -\Delta(\nabla^{\GI} E) &= -\grad\left(f\circ \varphi^{-1} (1-\sqrt{|D\varphi^{-1}|})\right) \nonumber \\
  & \hspace{10pt} - \sqrt{|D\varphi^{-1}|\,I_0\circ\varphi^{-1}} \grad \sqrt{I_1} + \grad\left(\sqrt{|D\varphi^{-1}|\,I_0\circ\varphi^{-1}}\right)\sqrt{I_1}
\end{align}
Since we are taking the Sobolev gradient of $E$, we apply the inverse Laplacian to the right hand side of the above equation to solve for $\nabla^{\GI} E$.
\end{proof}


\subsection{Thoracic Density Registration}\label{application_thoracic}
We now present application of the above developed theory to the problem of estimating complex anatomical deformations associated with the breathing cycle as imaged via Computed Tomography (CT)~\cite{gorbunova2012mass}. This problem has wide scale medical applications, in particular radiation therapy of the lung where accurate estimation of organ deformations during treatment impacts dose calculation and treatment decisions~\cite{sawant2014investigating,keall2005four,suh2014imrt,geneser2011quantifying}.
The current state-of-the-art radiation treatment planning involves the acquisition of a series of respiratory correlated CT (RCCT) images to build 4D (3 spatial and 1 temporal) treatment planning data sets.
Fundamental to the processing and clinical use of these 4D data sets is the accurate estimation of registration maps that characterize the motion of organs at risk as well as the target tumor volumes.

The 3D image produced from X-ray CT is an image of linear attenuation coefficients.
For narrow beam X-ray linear attenuation coefficient (LAC) for a single material (units $\mathrm{cm}^{-1}$) is defined as $\mu(x) = m\rho(x)$, where $m$ is a material-specific property called the mass attenuation coefficient (units $\mathrm{cm}^2/\mathrm{g}$) that depends on the energy of the X-ray beam. Linear attenuation coefficient is proportional to the true density and therefore exhibits conservation of mass. 
Unfortunately, CT image intensities do not represent true narrow beam linear attenuation coefficients.
Instead, modern CT scanners use wide beams that yield secondary photon effects at the detector.
CT image intensities reflect \emph{effective} linear attenuation coefficients as opposed to the true narrow beam linear attenuation coefficient.

\begin{figure}
  \centering                     
  \includegraphics[trim={1.0cm, 0.5cm, 1.8cm, 1.0cm},clip,width=.57\textwidth]{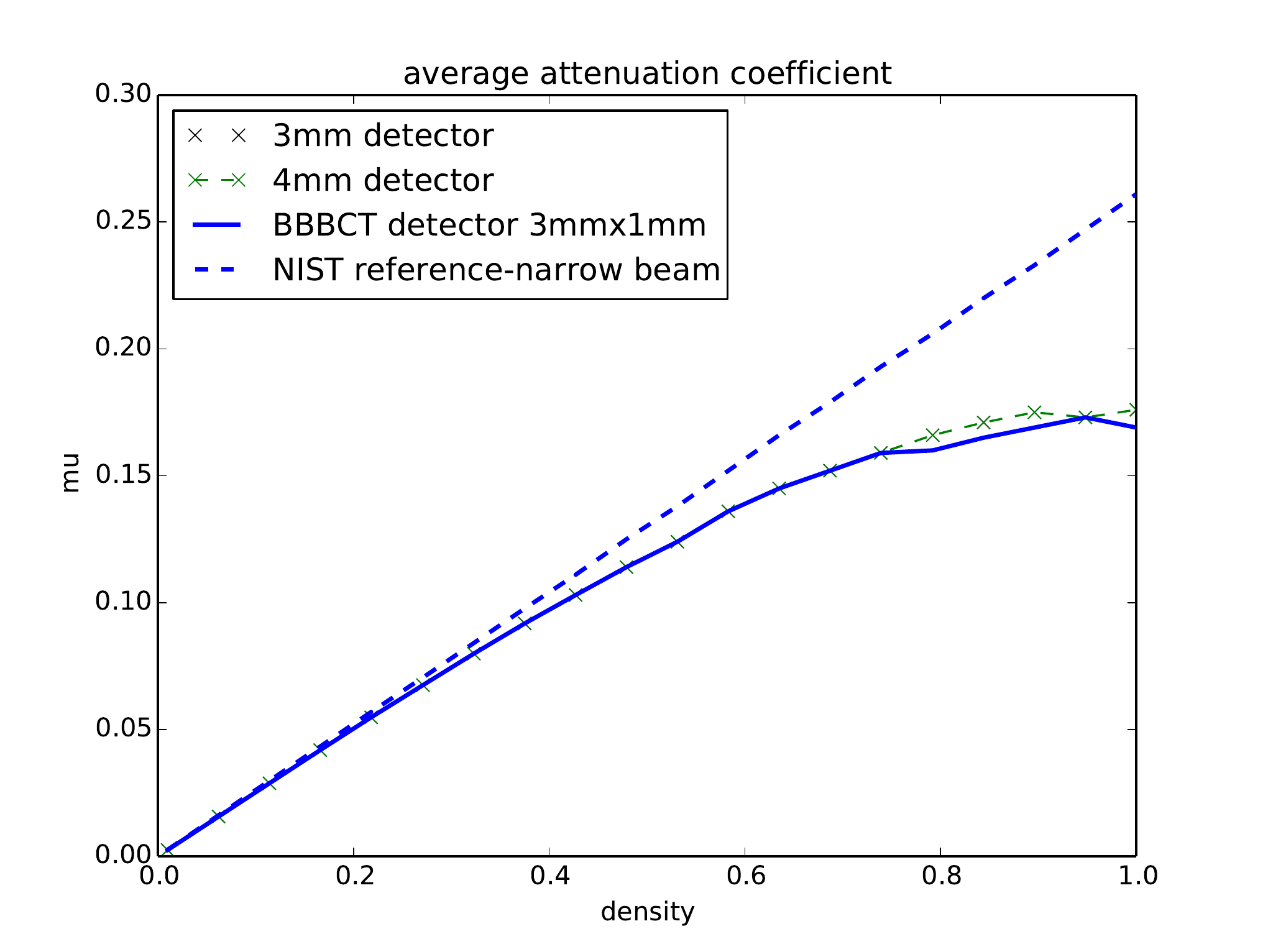}
  \caption{$\frac{1}{\alpha}$-Mass conservation in lung density matching: 
    Effective LAC from Monte Carlo simulation (solid line) and NIST reference narrow beam LAC (dashed line).
    The true relationship between effective LAC and narrow beam LAC is nonlinear.
  }
  \label{fig:montecarlo}
\end{figure}

To see the relationship between effective LAC and true narrow beam LAC, we ran a Monte Carlo simulation using an X-ray spectrum and geometry from a Philips CT scanner at various densities of water (since lung tissue is very similar to a mixture between water and air) \cite{boone1997accurate}.
The nonlinear relationship between effective LAC and narrow beam LAC relationship is clear (see Figure~\ref{fig:montecarlo}).

If we have conservation of mass within a single subject in a closed system, we expect an inverse relationship between average density in a region $\Omega$ and volume of that region: $D_t = \frac{M}{V_t}$.
Here $V_t = \int_{\Omega_t} 1 dx$, $D_t = \int_{\Omega_t} I_t(x) dx/V_t$, $\Omega_t$ is the domain of the closed system (that moves over time), and $t$ is a phase of the breathing cycle.
This relationship becomes linear in log space with a slope of $-1$:
\begin{equation}
  \ln(D_t) = \ln(M) - \ln(V_t)
\end{equation}

Our experimental results confirm the Monte Carlo simulation in that lungs imaged under CT do not follow this inverse relationship.
Rather, the slope found in these datasets in log space is consistently greater than -1 (see Figure~\ref{fig:regresults}).
This implies for real clinical CT data sets the lung tissue is acted on by an $\alpha-density$ action. Using the isomorphism  between $\alpha-densities$ and $1-densities$, we estimate a power transformation, i.e. $I(x) \mapsto I(x)^\alpha$, and estimate the $\alpha$ that yields the best conservation of mass property.

For each subject, we perform a linear regression of the measured LAC density in the homogeneous lung region and the calculated volume in log space. Let $\vect{d}(\alpha) = \log\left(\int_{\Omega_t} I_t(x)^\alpha dx/\int_{\Omega_t}1 dx \right)$ (the log density) and $\vec{v} = \log(\int_{\Omega_t}1 dx)$ (the log volume), where again $t$ is a breathing cycle timepoint.
The linear regression then models the relationship in log space as $\vect{d}(\alpha) \approx a \vec{v} + b$.
Let $a_j(\alpha)$ be the slope solved for in this linear regression for the $j^{th}$ subject.
To find the optimal $\alpha$ for the entire dataset, we solve
\begin{equation}
  \alpha = \argmin_{\alpha'} \sum_j (a_j(\alpha') + 1)^2,
\end{equation}
which finds the value of $\alpha$ that gives us an average slope closest to -1. We solve for $\alpha$ using a brute force search.

Applying this power function to the CT data allows us to perform our density registration algorithm based on the theory developed.
We therefore seek to minimize the energy functional described in the previous section, given by:
\begin{align}
  E(\varphi) &= d^{2}_F( \varphi_*(f\, dx), (f \circ \varphi^{-1}) dx) + d^{2}_F( \varphi_*(I_0\, dx), I_1\, dx)) \\
             &= \underbrace{\int_{\Omega} (\sqrt{|D\varphi^{-1}|}-1)^2\,f\circ\varphi^{-1} \,dx}_{E_1(\varphi)} +
               \underbrace{\int_{\Omega}\Big(\sqrt{|D\varphi^{-1}|I_0\circ\varphi^{-1}}- \sqrt{I_1} \Big)^2\,dx}_{E_2(\varphi)} \;.
\end{align}

We construct the density $f(x) \, dx$, a positive weighting on the domain $\Omega$, to model the physiology of the thorax: regions where $f(x)$ is high have a higher penalty on non-volume preserving deformations and regions where $f(x)$ is low have a lower penalty on non-volume preserving deformations.
Physiologically, we know the lungs are quite compressible as air enters and leaves.
Surrounding tissue including bones and soft tissue, on the other hand, is essentially incompressible.
Therefore, our penalty function $f(x)$ is low inside the lungs and outside the body and high elsewhere.
For our penalty function, we simply implement a sigmoid function of the original CT image: $f(x) = \mathrm{sig}(I_0(x))$.


Recalling the Sobolev gradient calculated in Theorem~\ref{Thm-H1Grad} with respect to the energy functional given by
\begin{align}
  \delta E &= -\Delta^{-1} \Big( -\nabla\big(f\circ \varphi^{-1} (1-\sqrt{|D\varphi^{-1}|})\big) - \nonumber \\
  & \hspace{30pt} \sqrt{|D\varphi^{-1}|\,I_0\circ\varphi^{-1}} \nabla\big(\sqrt{I_1}\big) + \nabla\big(\sqrt{|D\varphi^{-1}|\,I_0\circ\varphi^{-1}}\big)\sqrt{I_1} \Big)\;.
\end{align}
Then, the current estimate of $\varphi^{-1}$ is updated directly via a Euler integration of the gradient flow \cite{Rottman2016}:
\begin{equation}
  \varphi^{-1}_{j+1}(x) = \varphi_{j}^{-1}(x + \epsilon \delta E )
\end{equation}
for some step size $\epsilon$.
Since we take the Sobolev gradient the resulting deformation is guaranteed to be invertible with a sufficiently small $\epsilon$. Also notice that the gradient only depends on $\varphi^{-1}$ so there is no need to keep track of both $\varphi$ and $\varphi^{-1}$. The exact numerical algorithms is as follows:

\begin{tcolorbox}
\textbf{Numerical algorithm for Weighted Diffeomorphic Density Registration}

	\begin{algorithmic}
		\STATE Chose $\epsilon > 0$
	\STATE  $\varphi^{-1} \gets id$
	\STATE  $|D\varphi^{-1}| \gets 1$
	\FOR{$iter = 0 \cdots numiter $}
	\STATE  $\varphi_{*}I_{0} \gets  I_{0} \circ \varphi{-1} |D\varphi^{-1}| $
	\STATE $u \gets -\nabla \big(f\circ \varphi^{-1} (1-\sqrt{|D\varphi^{-1}|}) \big) - \sqrt{|D\varphi^{-1}|\,I_0\circ\varphi^{-1}} \nabla\big(\sqrt{I_1}\big)+ \nabla\big(\sqrt{|D\varphi^{-1}|\,I_0\circ\varphi^{-1}}\big)\sqrt{I_1} \Big)$
	\STATE $v \gets -\Delta^{-1}(u)$
	\STATE $\varphi^{-1}(y) \gets  \varphi^{-1}(y+\epsilon v)$
	\STATE $|D\varphi^{-1}| \gets |D\varphi^{-1}| \circ \varphi^{-1} e^{-\epsilon \mbox{div}(v)}$
	\ENDFOR
    \end{algorithmic}
\end{tcolorbox}
The algorithm was implemented using the  PyCA  package
and can be downloaded at
\begin{center}
	\texttt{\url{https://bitbucket.org/crottman/pycaapps/src/master/}}
\end{center}
See the application \emph{Weighted Diffeomorphic Density Registration}.

\begin{figure}[htbp]
  \centering
    \includegraphics[width=0.49\textwidth]{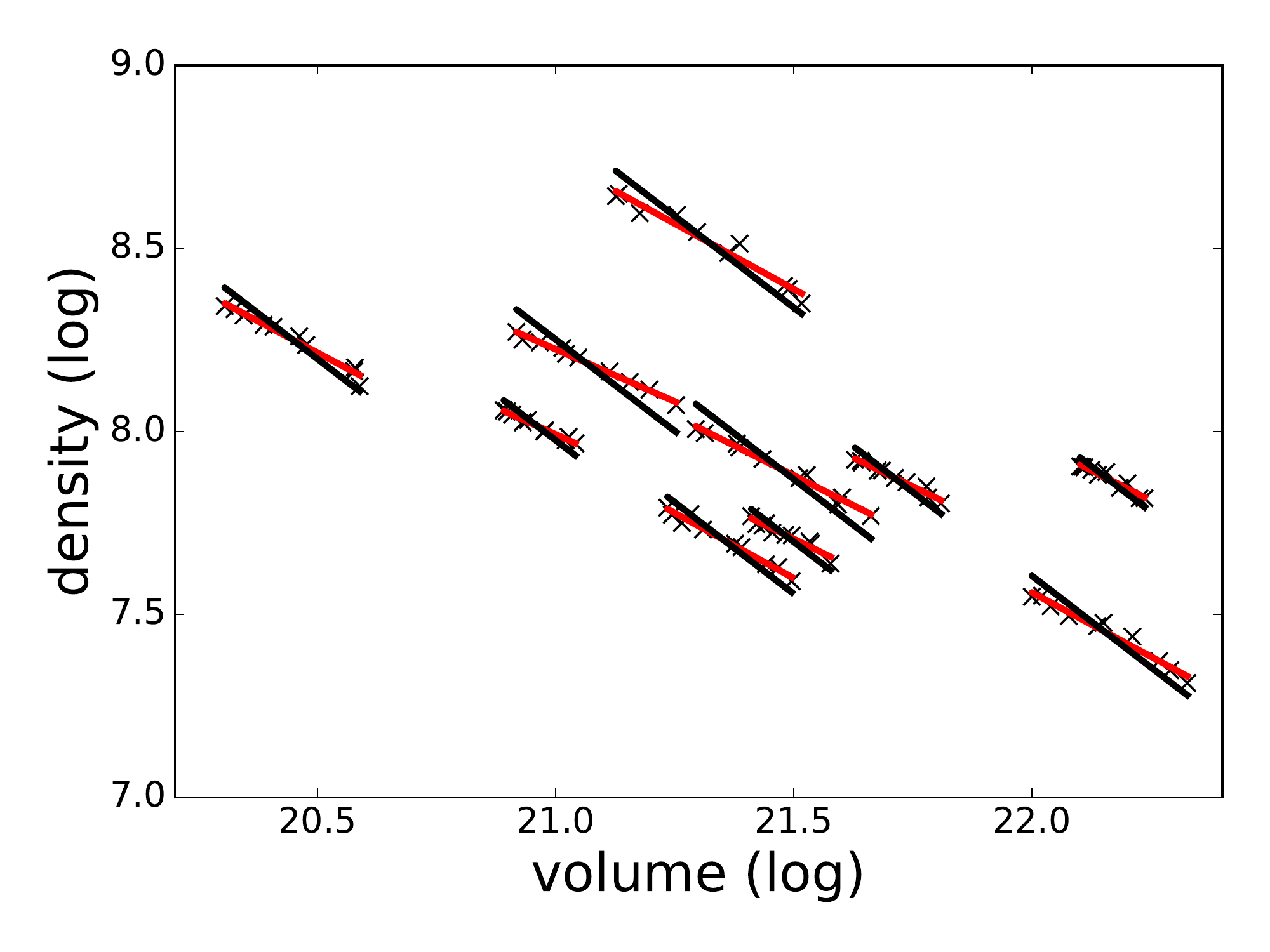}
    \includegraphics[width=0.49\textwidth]{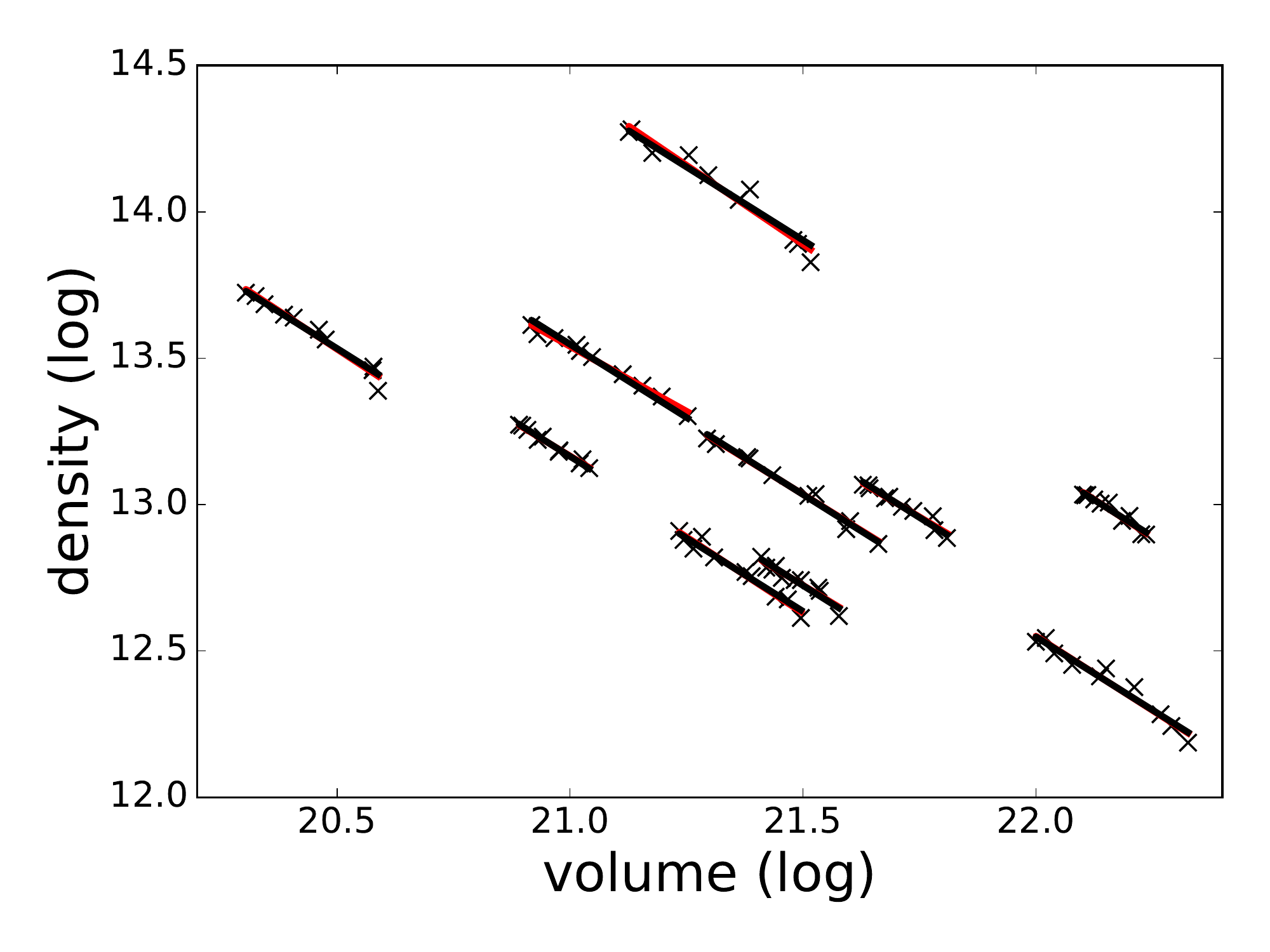}
    \includegraphics[trim={.5cm, 0.1cm, .5cm, 7cm},clip,width=0.49\textwidth]{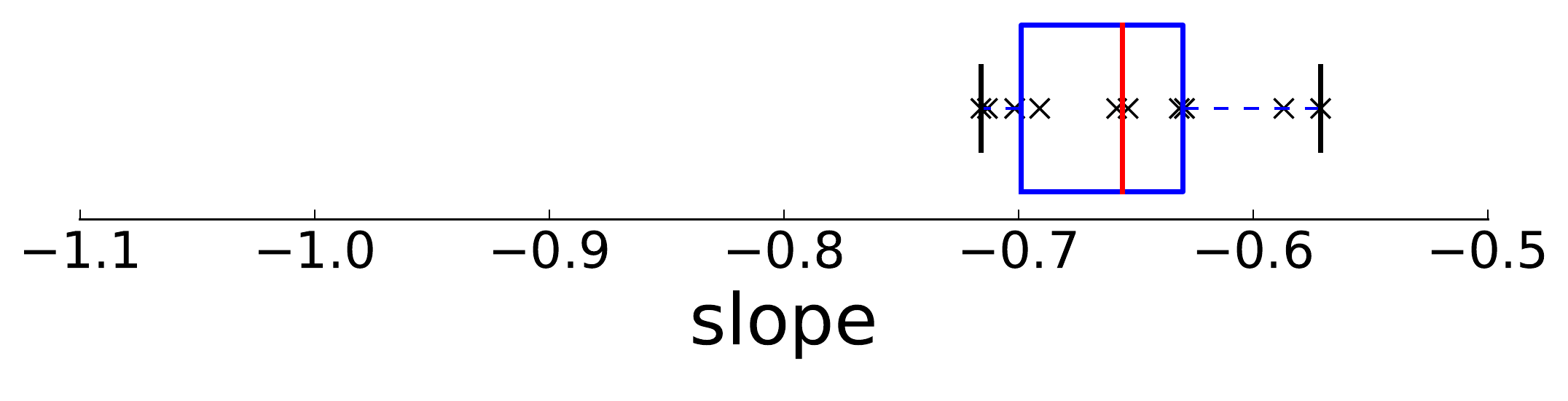}
    \includegraphics[trim={.5cm, 0.1cm, .5cm, 7cm},clip,width=0.49\textwidth]{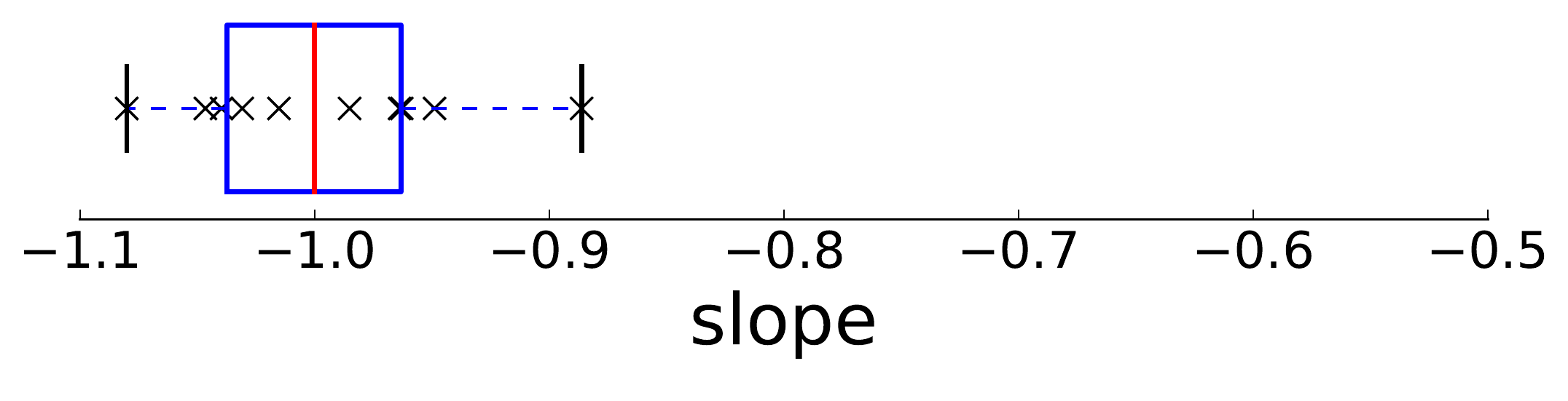}
\caption{Density and volume log-log plots.
  Upper left: log-log plots without applying the exponential correction for all ten DIR subjects.
  The best fit line to each dataset is in red and the mass-preserving line (slope = -1) is in black.
  Upper right: log-log plots after applying the exponential correction $I(x)^\alpha$ to the CT images.
  In this plot, the best fit line matches very closely to the mass-preserving line.
  Bottom row: corresponding box plots of the slopes found in the regression.
  }
\label{fig:loglog}
\end{figure}

\begin{figure}[htbp]
  \centering                    
    \includegraphics[trim={4.3cm, 2cm, 4.2cm, .5cm},clip,width=0.31\textwidth]{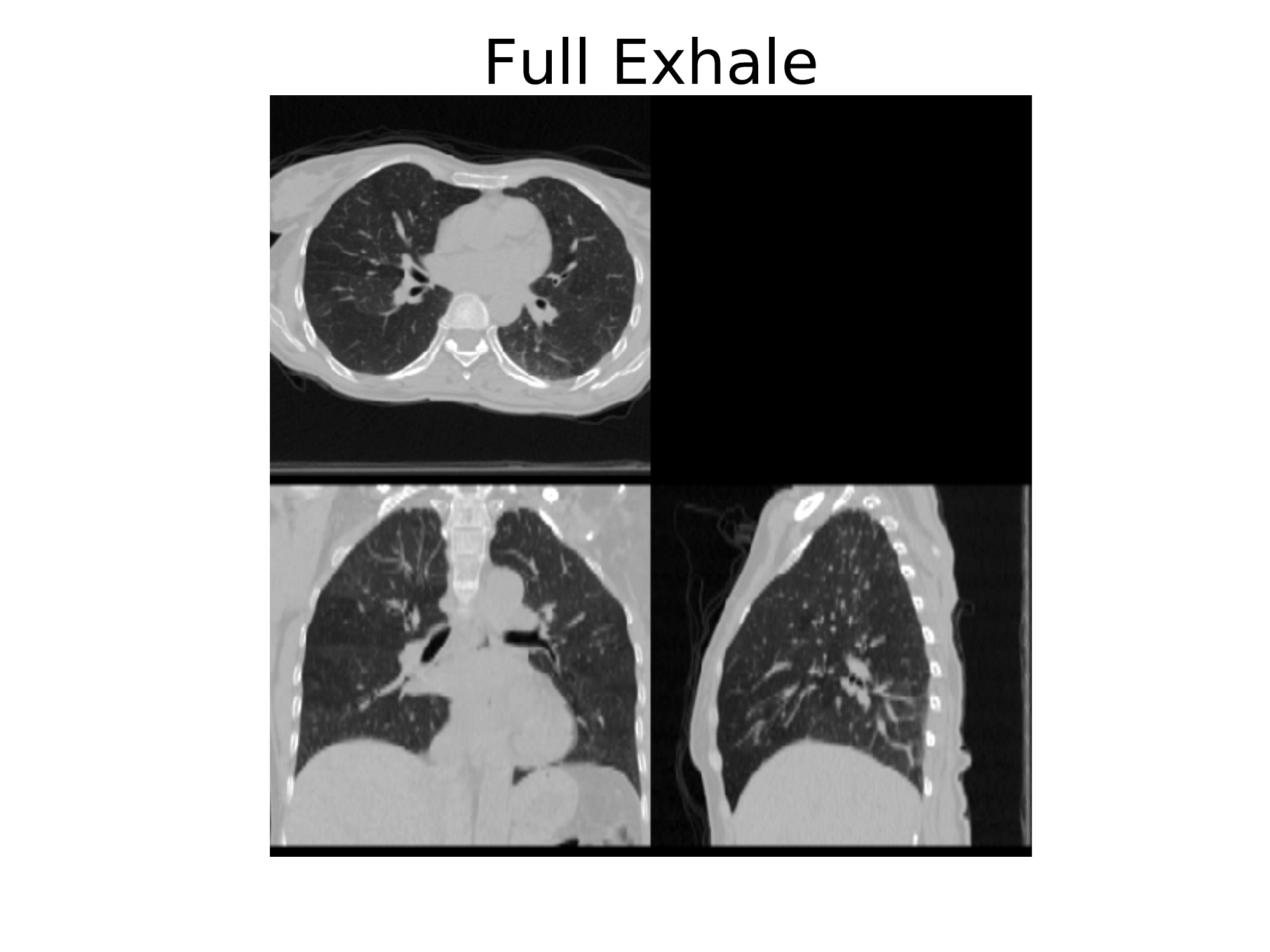}
    \includegraphics[trim={4.3cm, 2cm, 4.2cm, .5cm},clip,width=0.31\textwidth]{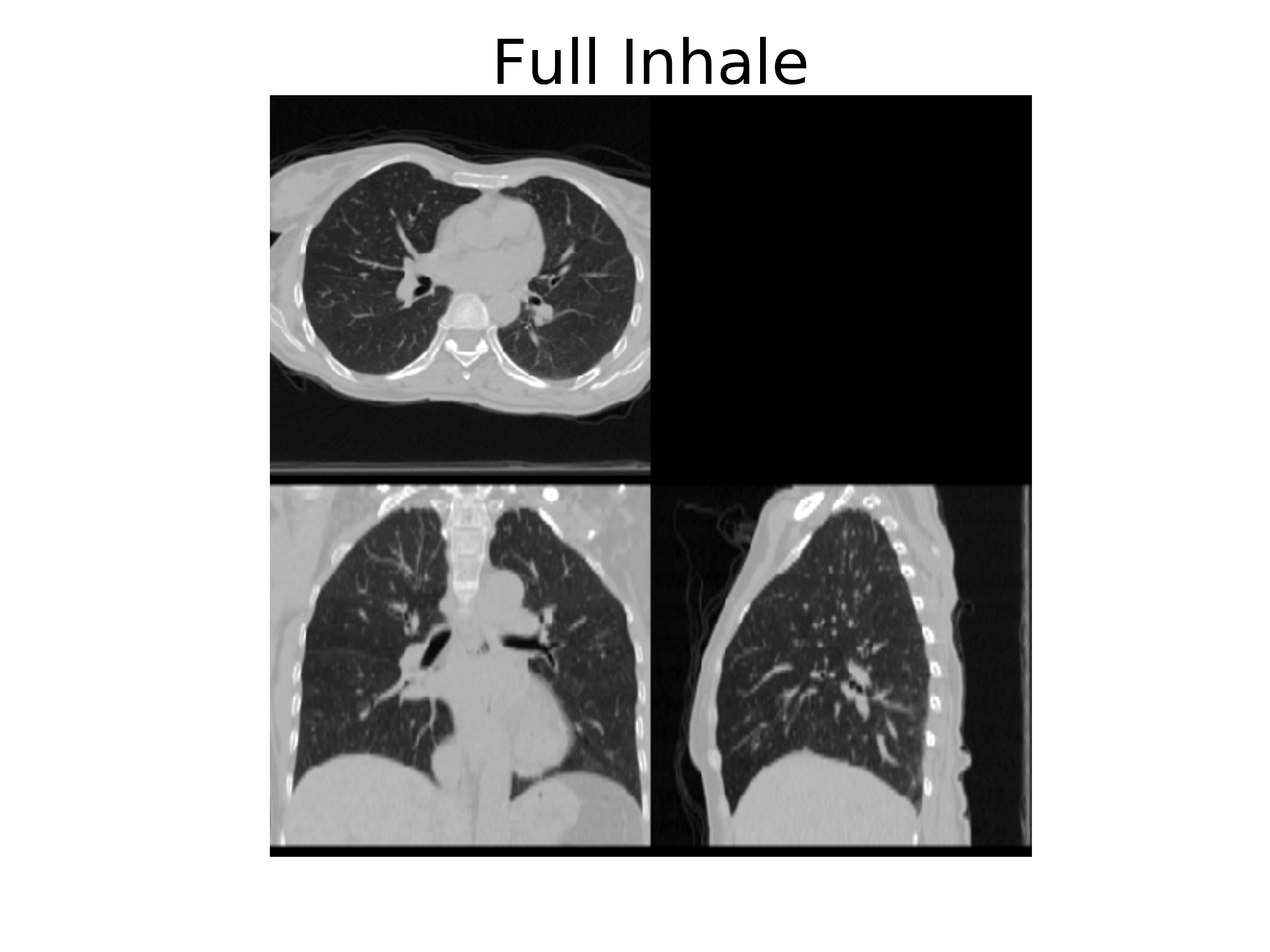}
    \includegraphics[trim={4.3cm, 2cm, 4.2cm, .5cm},clip,width=0.31\textwidth]{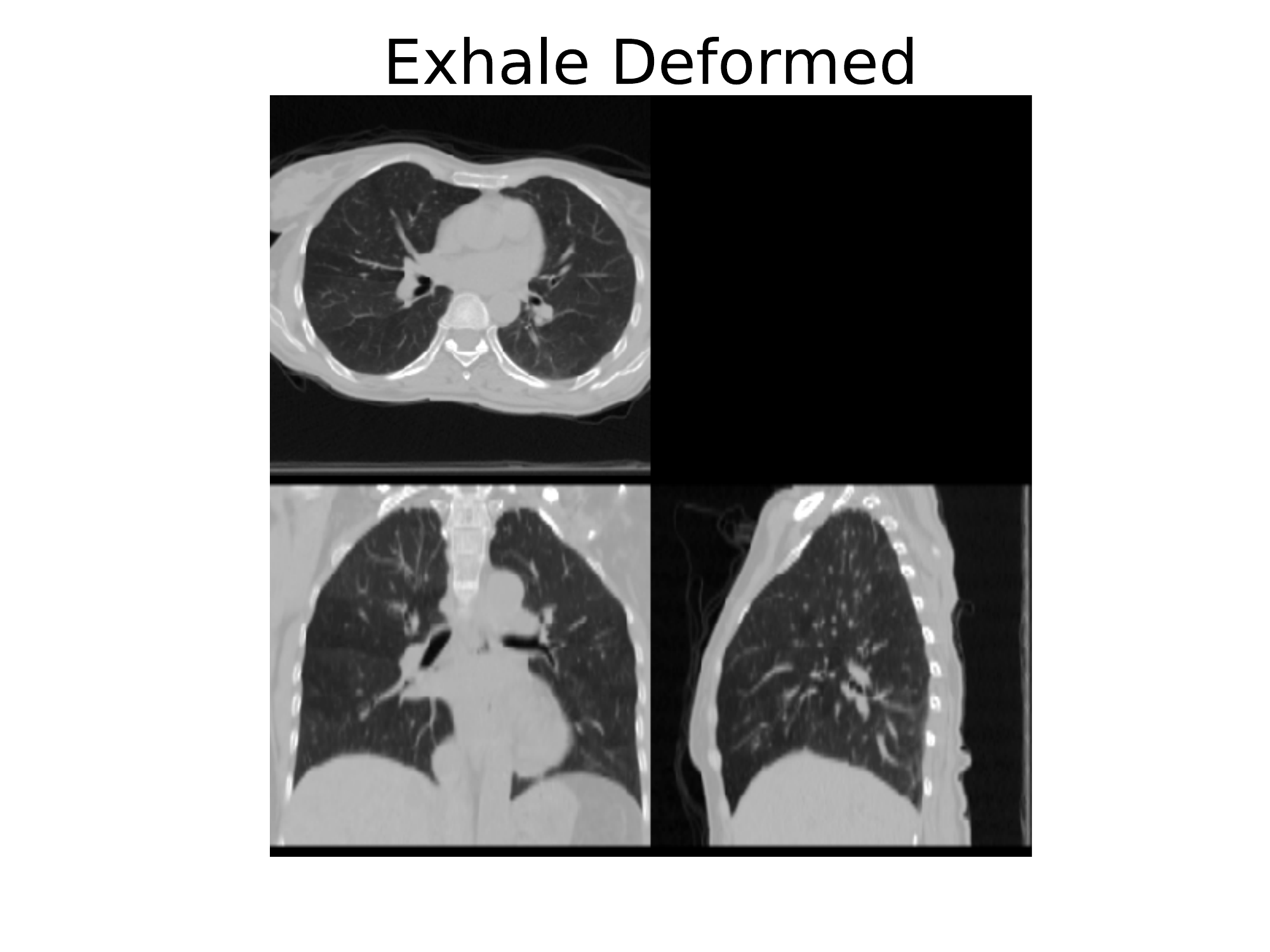}
    \includegraphics[trim={3cm, 1.0cm, 2.0cm, .3cm},clip,width=0.31\textwidth]{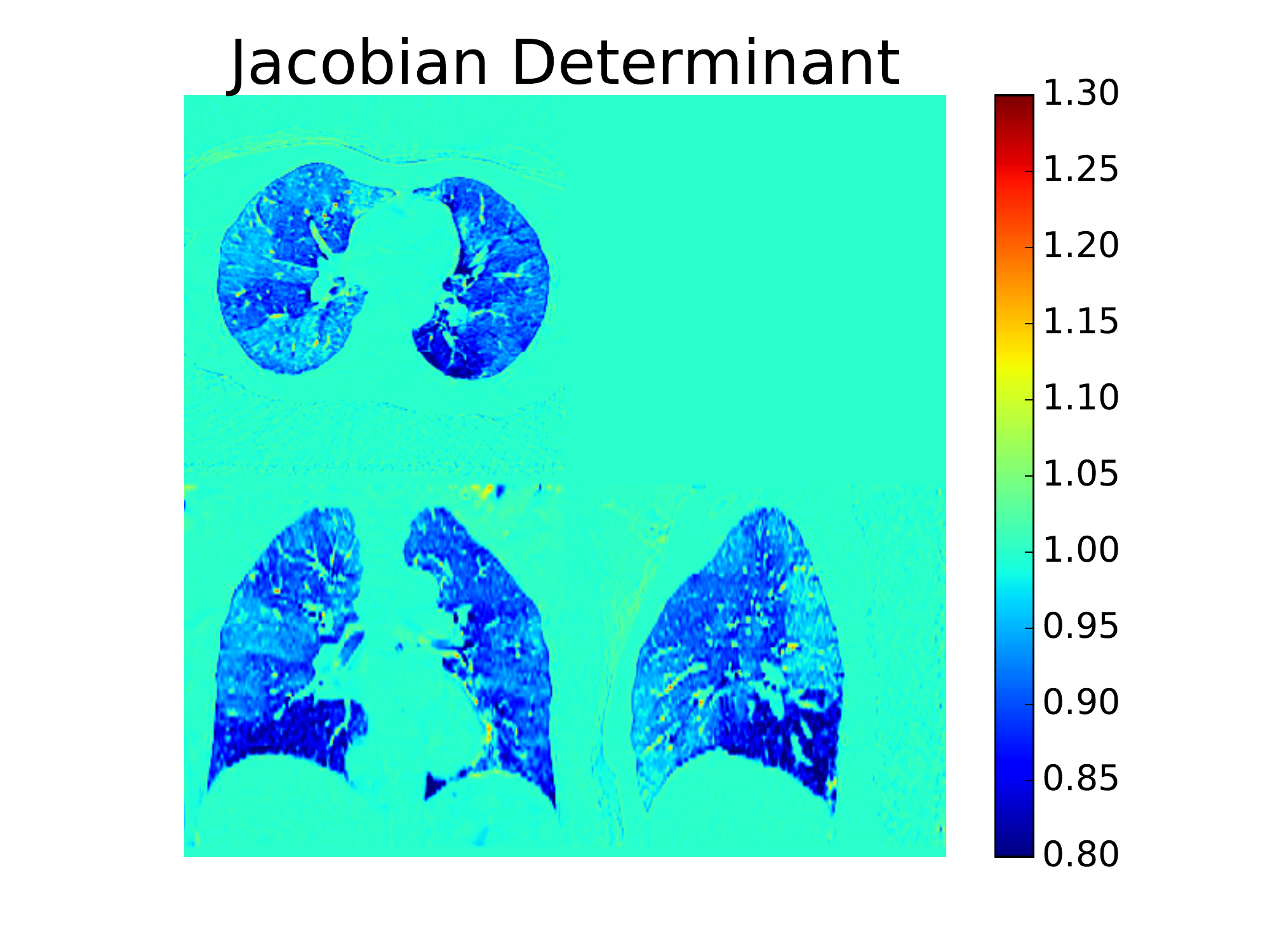}
    \includegraphics[trim={3cm, 1.0cm, 2.0cm, .3cm},clip,width=0.31\textwidth]{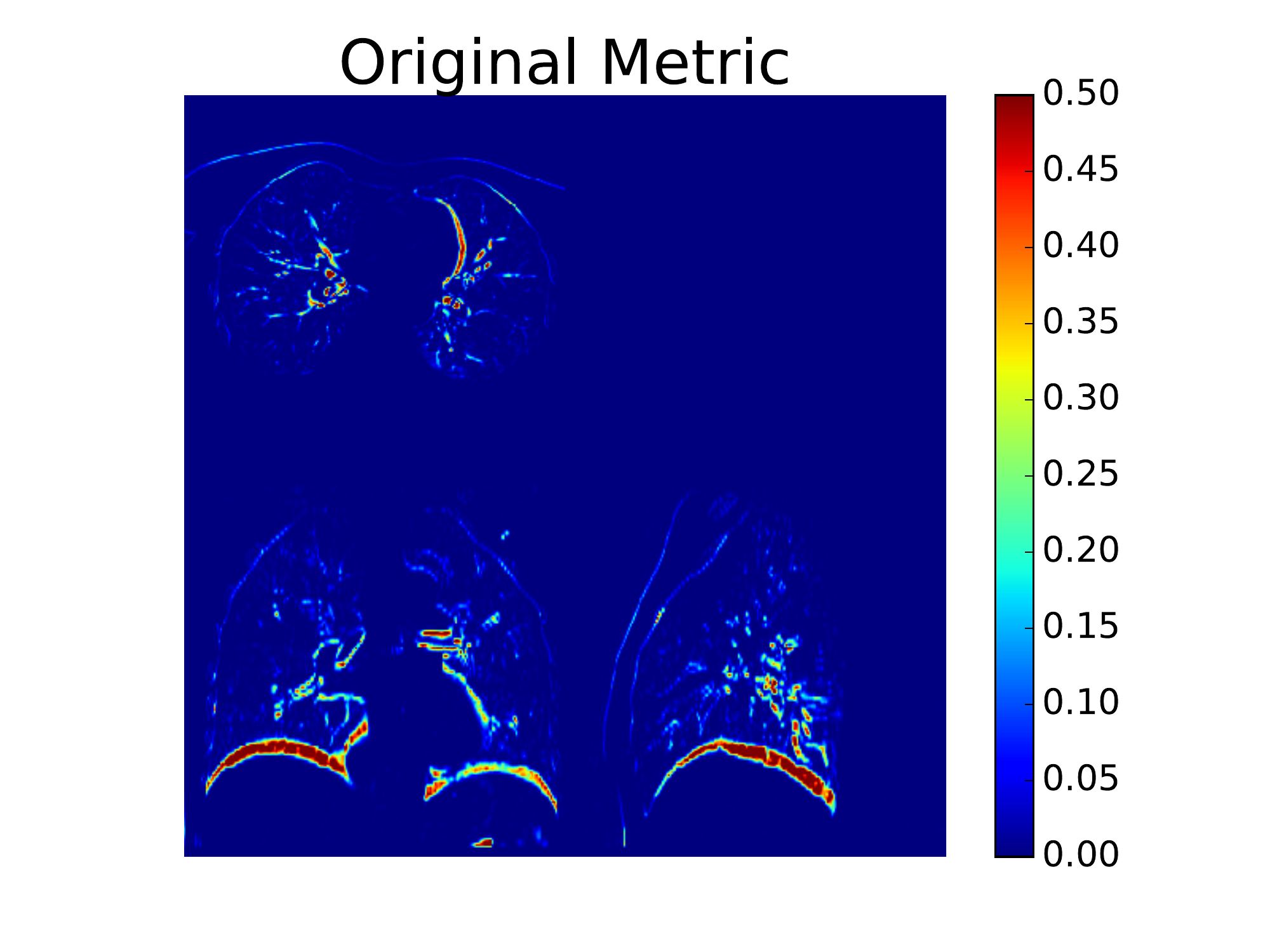}
    \includegraphics[trim={3cm, 1.0cm, 2.0cm, .3cm},clip,width=0.31\textwidth]{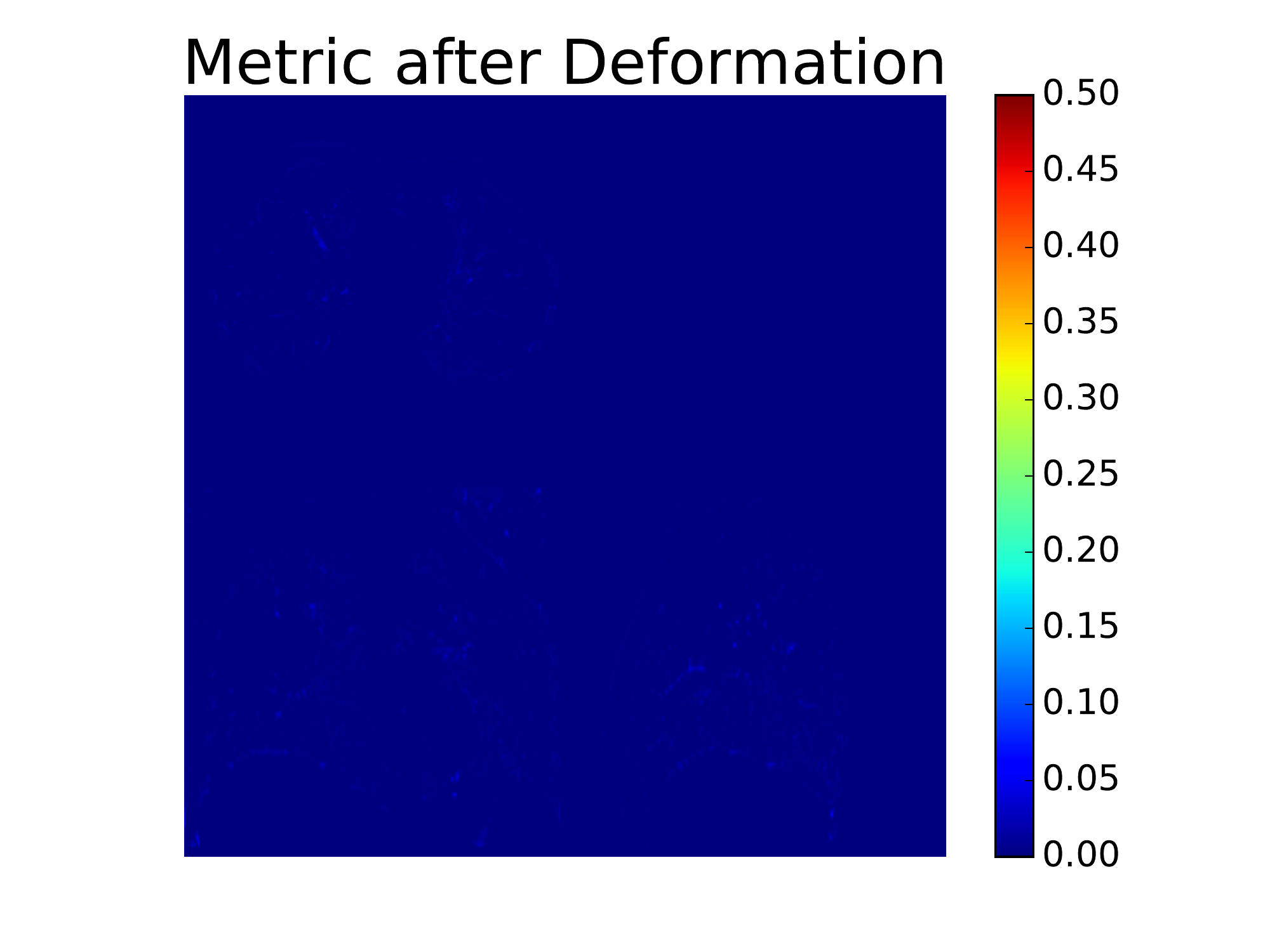}
    \includegraphics[trim={0.5cm, 1.0cm, 1.0cm, .5cm},clip,width=0.38\textwidth]{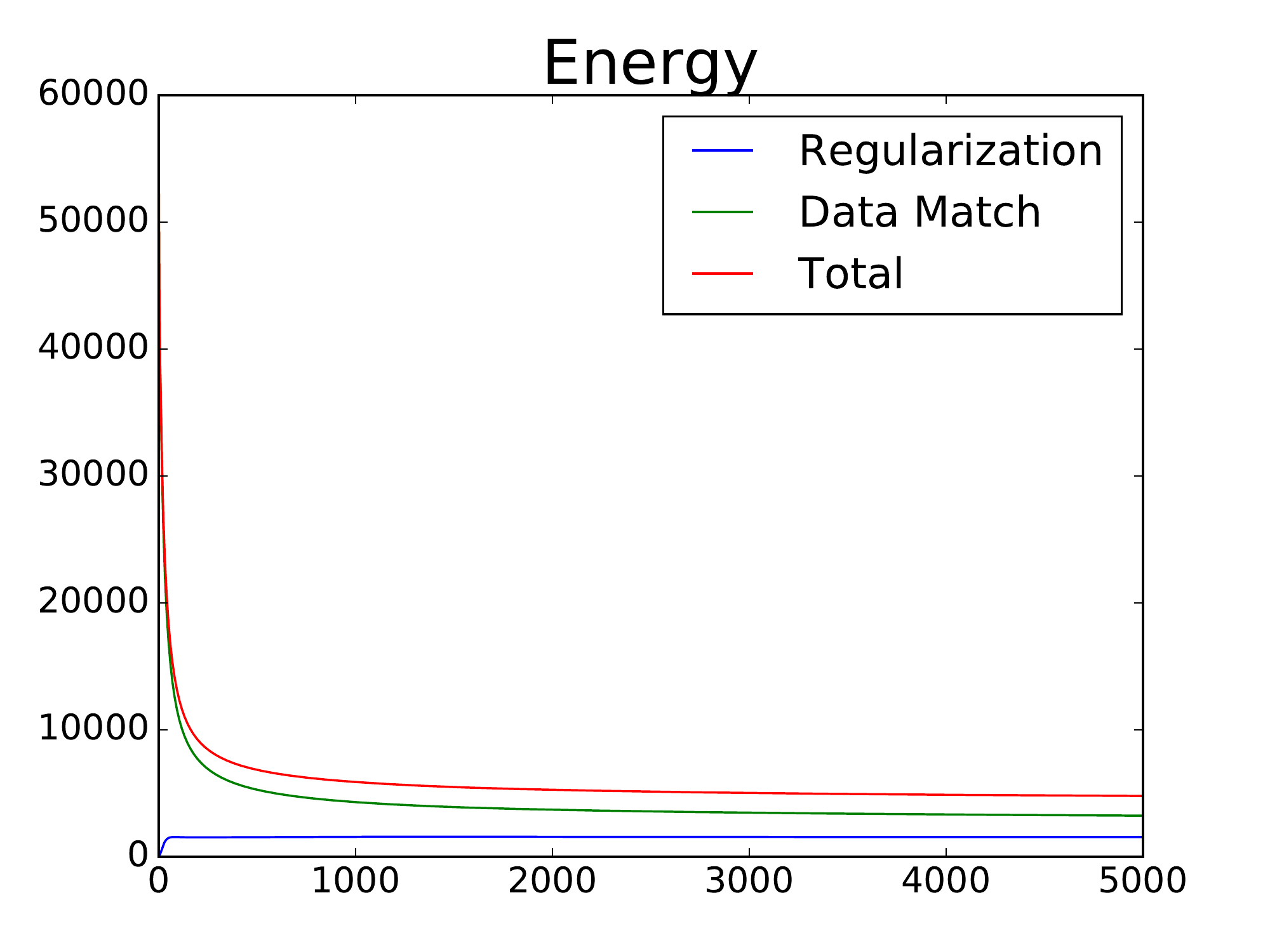}
    \includegraphics[trim={3cm, 1.0cm, 2.0cm, .3cm},clip,width=0.31\textwidth]{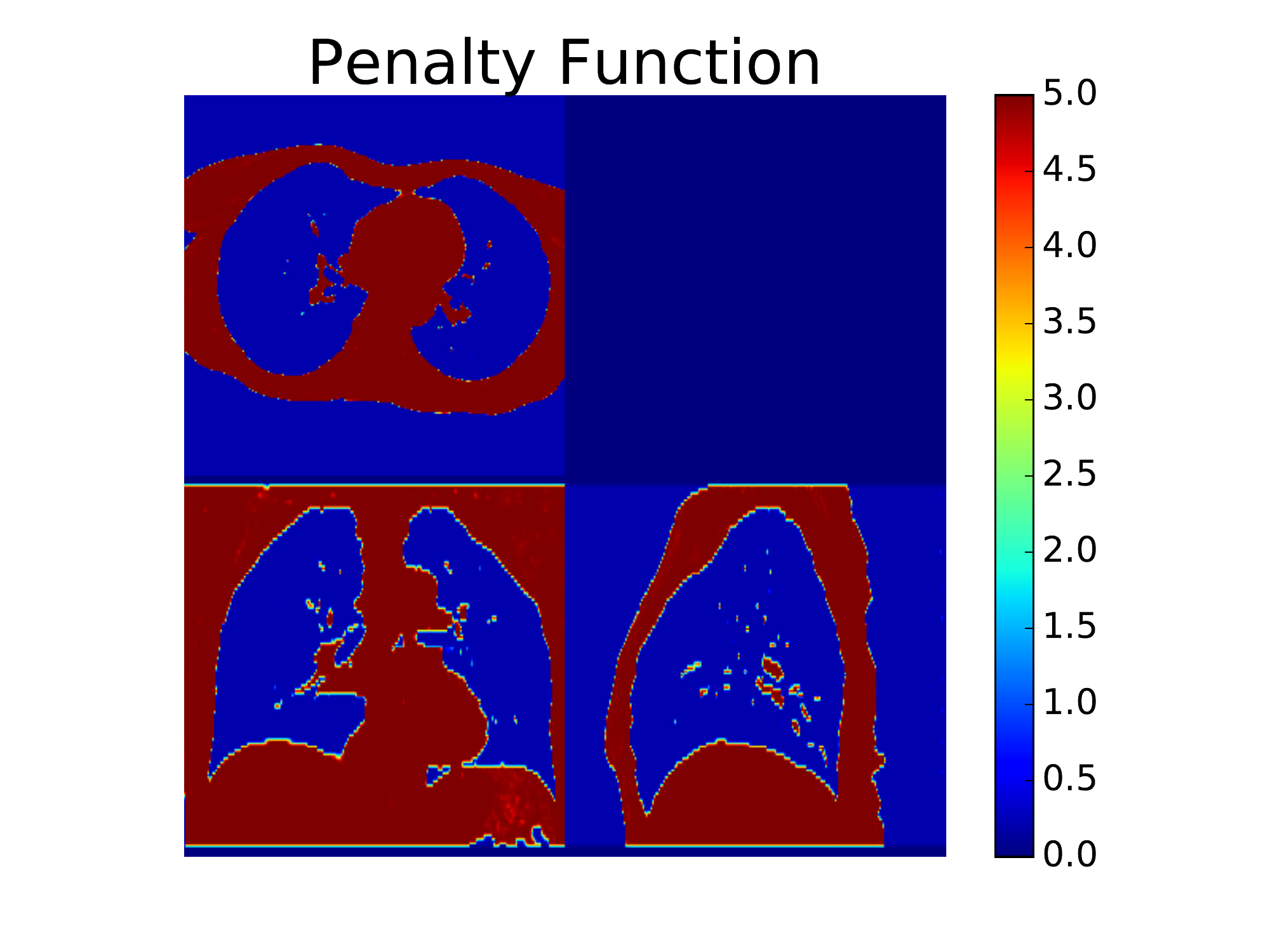}
\caption{Registration results.
  Top row: full inhale, full exhale, and the deformed exhale density estimated using our method.
  Middle row: Jacobian determinant of the transformation, initial Fisher-Rao metric, and Fisher-Rao metric after applying the density action.
Notice that outside the lungs the estimated deformation is volume preserving.
  Bottom row: Energy as a function of iterations, and penalty function.
  }
\label{fig:regresults}
\end{figure}

For the DIR dataset, we solved for the exponent that yields conservation of mass, which yielded  
 $\alpha=0.60$ that gives us the best fit.
Without using the exponential fit, the average slope of log density log volume plot was -0.66 (SD 0.048).
After applying the exponential to the CT intensities, the average slope is -1.0 (SD 0.054).
The log-log plots of all ten patients in the DIR dataset as well as box plots of the slope is shown in Figure~\ref{fig:loglog}.

For the 30 subject dataset, we solved for $\alpha = 0.52$ that gives us conservation of mass.
Without using the exponential fit, the average slope of the log-log plot was -0.59 (SD 0.11).

We applied our proposed weighted density registration algorithm to the first subject from the DIR dataset.
This subject has images at 10 timepoints and has a set of 300 corresponding landmarks between the full inhale image and the full exhale image.
These landmarks were manually chosen by three independent observers.
Without any deformation, the landmark error is 4.01 mm (SD 2.91 mm).
Using our method, the landmark error is reduced to 0.88 mm (SD 0.94 mm), which is only slightly higher than the observer repeat registration error of 0.85 mm (SD 1.24 mm).

We implement our algorithm on the GPU and plot the energy as well as the Fisher-Rao metric with and without applying the deformation.
These results are shown in Figure~\ref{fig:regresults}.
In this figure, we show that we have excellent data match, while the deformation remains physiologically realistic: inside the lungs there is substantial volume change due to respiration, but the deformation outside the lungs is volume preserving.
With a $256 \times 256 \times 94$ voxel dataset, our algorithm takes approximately nine minutes running for four thousand iterations on a single nVidia Titan Z GPU.
\subsection*{Acknowledgements}
We thank Caleb Rottmann, who worked on the implementation of the weighted diffeomorphic density matching algorithm. We are grateful for valuable discussions with Boris Khesin, Peter Michor  and FX Vialard. 

This work was partially supported by the grant NIH R01 CA169102, the Swedish
Foundation for Strategic Research (ICA12-0052), an EU Horizon 2020 Marie
Sklodowska-Curie Individual Fellowship (661482) and by the Erwin Schr\"odinger
Institute programme: Infinite-Dimensional Riemannian Geometry with Applications
to Image Matching and Shape Analysis by the FWF-project P24625.

\end{document}